\documentclass[10pt,english,a4paper]{article}
\usepackage{amssymb}
\usepackage{amsmath}
\usepackage{geometry}
\usepackage{graphics}
\usepackage{epsfig}
\usepackage{url}
\usepackage{array}
\usepackage{amsthm}
\usepackage{chicago} 
\usepackage[normalem]{ulem}
\usepackage{color}
\usepackage{multirow}
\definecolor{mygrey}{gray}{0.35}

\def\eps{\varepsilon}
\def\chi{v}

\def\mygrey#1{\color{mygrey}#1\color{black}}

\newtheorem{theorem}{Theorem}
\newtheorem{definition}{Definition}
\newtheorem{lemma}{Lemma}

\newcommand{\highlight}[1]{#1}

\begin{document}
\date{\today}

\title{Heuristic and exact solutions to the inverse power index problem for small voting bodies$^\dagger$}
\author{Sascha Kurz and Stefan Napel\\
        University of Bayreuth, 95440 Bayreuth, Germany\\
        \{sascha.kurz,stefan.napel\}@uni-bayreuth.de}

\begin{abstract}
Power indices are mappings that quantify the influence of the members of a voting body on collective decisions a~priori. Their nonlinearity
and discontinuity makes it difficult to compute inverse images, i.e., to determine a voting system which induces a power distribution as
close as possible to a desired one. \highlight{The paper considers approximations to this inverse problem for the Penrose-Banzhaf
index by hill-climbing algorithms and exact solutions which are obtained by enumeration and integer linear programming techniques.} They are compared to the results of three simple solution
heuristics. The heuristics perform well in absolute terms but can be improved upon very considerably in relative terms. The findings complement
known asymptotic results for large voting bodies and may improve termination criteria for local search algorithms. 

  \bigskip

 \noindent
 \textbf{Keywords} electoral systems; simple games; weighted voting games; square root rule; Penrose limit theorem; Penrose-Banzhaf index; institutional design
 
  \noindent
  \textbf{Mathematics Subject Classification (2010)} 91B12, 91A12, 90C10
\end{abstract}


\maketitle

\section{Introduction}\label{sec:intro}

Collective decision rules and, in particular, heterogeneous voting weights for members of a committee, council, or shareholder meeting
translate into \emph{influence} on collective decisions in a nonlinear and even discontinuous fashion. This can be seen, for instance,
by considering a decision quota of $q=50\%$ and players $i\in N=\{1,2,3\}$ whose voting weights are given by either the vector
(i) $w=(w_1,w_2,w_3)=(33.\bar 3, 33.\bar 3, 33.\bar 3)$, (ii) $w'=(50-\eps, 48+\eps, 2)$, or (iii) $w''=(50+\eps, 48-\eps, 2)$ for
small $\eps>0$. The major weight change from $w$ to $w'$ does not affect possibilities to form a \emph{winning coalition} at all,
where coalition $S\subseteq N$ is called winning if 
the cumulative weight of its members exceeds the quota. Namely, $S$ is winning if and only if $|S|\geq 2$. By symmetry,
the distribution of influence can a~priori be expected to equal $(\frac{1}{3},\frac{1}{3},\frac{1}{3})$ for either of the voting
systems described by $(q;w)$ and $(q;w')$. The minor change from $w'$ to $w''$, in contrast, renders voter~1 a dictator with
associated power distribution $(1,0,0)$.

Social scientists, philosophers and mathematicians have investigated various \emph{voting power indices} which try to quantify
the a~priori distribution of influence on committee decisions. The Shapley-Shubik index \cite{Shapley/Shubik:1954} and the
Penrose-Banzhaf index (PBI) \cite{Penrose:1946,Banzhaf:1965} are most prominent, but by far not the only ones.\footnote{See
\citeN{Felsenthal/Machover:1998} or \citeN{Laruelle/Valenciano:2008} for overviews.} They help researchers clarify the non-trivial
a~priori power implications of different voting weight assignments to a wider audience. The combinatorial nature of weighted voting
systems can \highlight{easily} mislead the general public's intuition and even that of political practitioners. For instance, it was apparently
not noted that the voting weights of the original six members of the European Economic Community, in use from 1958 to 1973, rendered
Luxembourg a \emph{null player} whenever the EEC Council applied its qualified majority rule, \highlight{i.e., the country's vote was
irrelevant for the collective decision}. The public discussion~-- very heated in, but not restricted to, Poland and Germany~-- in the
wake of the 2007 EU~summit at which new voting rules for the EU~Council were agreed reflected persistent confusion between voting
weights and power.

Even to specialists, the discrete nature of voting rules still poses challenges. This is true, in particular, for the optimal design
of a voting system. Certain normative desiderata, such as the equal representation of bottom-tier voters in a two-tier voting system,
often call for a specific distribution of voting power. It is then a non-trivial exercise to find a deterministic voting rule that
comes as close as possible to inducing the desired a~priori power distribution.\footnote{Non-deterministic rules such as
\emph{random dictatorship} or also \emph{random quota rules} \cite[sec.~5]{Dubey/Shapley:1979} can easily solve the problem, but are
generally not regarded as satisfactory.} \highlight{Simple 
hill-climbing algorithms, such as the ones considered by \citeN{Laruelle/Widgren:1998},} \highlight{Leech \citeyear{Leech:2002,Leech:2002:IMF,Leech:2003},}
\citeN{Aziz/Paterson/Leech:2007}, or \citeN{Fatima/Wooldridge/Jennings:2008}, deliver excellent results for many instances
of this so-called \emph{inverse power index problem} but have \highlight{rarely been evaluated in a systematic fashion}.\footnote{\highlight{We thank an anonymous
referee for pointing out to us that \citeN[sec.~6]{Laruelle/Widgren:1998} contains the earliest published specification of an
iterative algorithm for solving the inverse problem. This algorithm has recently been evaluated and improved by \citeN{deNijs/Wilmer:2012}. \citeN[p.~44]{Lucas:1992} is one of the first to describe the basic approach
to implementing a target PBI vector by a suitable voting rule; while \citeN[p.~206f]{Nurmi:1982} seems to be the first to state
this institutional design problem and to discuss its relevance.}} One can neither rule out that only a
\emph{local minimum} of the distance between the desired and the induced power vector has been identified. Nor are \emph{bounds} known
on the possible gap to a globally optimal voting rule. The latter might involve the intersection of several one-dimensional $(q;w)$-rules
and, therefore, need not even be a feasible result of the applied search algorithm.

\highlight{For very big $n$ the distinctions between voting weight and
voting power tend to become negligible under simple majority rule if each voter's relative weight is small and vanishing in the limit (and one stays outside a class of somewhat pathological examples).\footnote{\highlight{The case in which the relative weights
of some voters do \emph{not} vanish as $n\to \infty$~-- so that a few large voters or \emph{atomic players} stick out in an ``ocean''
of infinitesimal voters (defining a so-called \emph{oceanic game})~-- reduces to the analysis of a modified voting body $v'$ for the PBI. It involves
only the atomic voters and a quota $q'$ which is obtained from $q$ by subtracting half of the aggregate weight of the infinitesimal voters
(see \citeNP[sec.~9]{Dubey/Shapley:1979}). If the number of atomic voters is small then the corresponding inverse problem is exactly
of the type investigated here.}\label{fn:oceanicgames} } \label{rem:pathological} Corresponding asymptotic results
for $n\to \infty$ which render the inverse problem trivial have first been suggested in the work of Lionel S.\ Penrose
\citeyear{Penrose:1946,Penrose:1952}. Rigorous investigations of asymptotic proportionality of voting weight and power, which is
often referred to as the \emph{Penrose limit theorem}, have only more recently been provided by \citeN{Neyman:1982}, \citeN{Lindner/Machover:2004}, \citeN{Chang/Chua/Machover:2006}, and \citeN{Lindner/Owen:2007}.}%

Motivated by qualified majority voting in the EU, S\l{}omczy\'nski and \.Zyczkowski
\citeyear{Slomczynski/Zyczkowski:2006,Slomczynski/Zyczkowski:2007} have identified an elegant way to approximately solve the inverse
problem \highlight{for $n$ still very distant from infinity} if the decision quota $q$ is a free parameter \highlight{and all
individual voting power targets are small compared to $\sqrt{n}$ times their quadratic mean}.
Their heuristic suggestion is particularly appealing because it avoids discrepancies between voting weight and power. Namely, approximate
proportionality between the normalized weight vector $w=(w_1, \ldots, w_n)$ and the induced PBI $B(q;w)$ is achieved \highlight{even for
$n\ll \infty$} when the quota is set to $q^*=\frac{1}{2}(1+\sqrt{\sum w_i^2})$. A desired power vector $\beta$ can hence approximately be
induced simply by choosing $w=\beta$ and then calculating $q^*$. Because the rule $(q^*;w)$ is simple and minimizes \highlight{unwanted} confusion
between voting weight and power, it has motivated the prominent ``Jagiellonian Compromise'' (also known as \emph{double square root voting
system}) in the discussion of future voting rules for the EU Council (see, e.g., \citeNP{Kirsch/Slomczynski/Zyczkowski:2007}).

Whether the decision quota $q$ is a free parameter, so that S\l{}omczy\'nski and \.Zyczkowski's optimal $q^*$ indeed can be chosen,
or not, depends on the application at hand. Even if it can, the lack of bounds on how well the $(q^*;\beta)$-heuristic performs relative
to the respective globally optimal solution to the inverse problem provides motivation for further research. Knowing more about the quality
of the $(q^*;\beta)$-heuristic is especially important for situations in which the heuristic can be expected to perform rather
badly. Specifically, the derivation of $q^*$ is based on a continuous approximation of the fundamentally discrete distribution of the
cumulative weight of a random coalition. Its use is problematic when this approximation is inaccurate. This pertains particularly~-- but
not exclusively~-- to ``small'' voting bodies with few members.

For a given number $n$ of players, the set of different binary voting systems or \emph{simple games} is finite. This finiteness entails
that many desired power distributions cannot be approximated too well. Nontrivially, this remains
true even for large $n$: \highlight{recent results by} Alon and Edelman \citeyear{Alon/Edelman:2010} imply \highlight{the existence of
sequences} of desired power distributions $\{\beta^n\}_{n=1, 2, \ldots}$ which stay at least a constant positive distance away (in the
$\Vert\cdot\Vert_1$-norm, \highlight{which adds up deviations from target for all voters $i=1, \ldots, n$)} from \emph{any} Penrose-Banzhaf power
distribution.

The finiteness of the set of  simple games \highlight{at the same time} suggests a trivial algorithm for
solving the inverse problem: enumerate all systems $\chi$ with $n$ players, compute the respective power distribution~-- say, the PBI
$B(\chi)$~-- and then pick a game $\chi^*$ that induces the smallest achievable difference between ideal vector $\beta$ and $B(\chi)$
according to a suitable \highlight{measure of distance (e.g., the metric induced by a particular vector norm)}.

To this end, a growing literature has investigated methods for the efficient enumeration of voting systems (see, e.g.,
\citeNP{Keijzer:2009}; \citeNP{Keijzer/Klos/Zhang:2010}; \citeNP{Kurz:2012:representation}). But, up to now, even the \emph{number} of
complete simple games (and also of weighted voting games) is unknown for $n>9$. So enumeration works only for voting bodies with very few members.
Exact solutions to the inverse problem can, fortunately, also be obtained for somewhat larger $n$ by \emph{integer linear
programming} (ILP) techniques. Such an approach was recently presented in \citeN{Kurz:2012}. It stands in the tradition of earlier
applications of ILP to electoral systems, as discussed in \citeN{Grilli/Manzi/Pennisi/Ricca/Simeone:1999},
\citeN{Pennisi/Ricca/Serafini/Simeone:2007} or \citeN{Ricca/Scozzari/Serafini/Simeone:2012}.

\highlight{The key idea is to translate the problem of finding a game $v$ whose PBI vector has a distance no greater than a given value
$\alpha\ge 0$ from the target vector $\beta$ into a set of linear inequality constraints plus the integer requirement that each coalition
is either winning or losing, and then to use efficient ILP solver software in order to check if a solution exists. A bisection process on $\alpha$, which consecutively halves the range of tentative optimal values of $\alpha$ that have not yet been ruled out by the computations, can be stopped after finitely many iterations 
(namely, once $|\alpha_{t+1}-\alpha_t|$ has fallen below the theoretical minimum distance between any two distinct $n$-voter PBI vectors).}

This paper draws on complete enumeration, \highlight{the indicated ILP bisection method, and standard hill-climbing} algorithms in order to evaluate the accuracy of three heuristic solutions
to the inverse problem for the PBI. The first heuristic simply combines $w=\beta$ with $q^{\circ}=50\%$; the second combines it with
the ``optimal quota'' $q^*$ derived by  \citeN{Slomczynski/Zyczkowski:2007}; the third uses $\bar q=\frac{1}{2}+{1}/{\sqrt{\pi n}}$.
The latter quota is the average of $q^*$ computed over a set of $\beta$-vectors which is of particular interest for the egalitarian
design of two-tier voting systems \cite{Slomczynski/Zyczkowski:2011}.

We compute differences between the \highlight{respective heuristic and exact solutions} for three different metrics and a comprehensive grid
of \highlight{conceivable} target vectors with up to $n=7$ voters. We study rather time-consuming approximations of the exact solutions for a large
sample of grid points for \highlight{$8\le n\le 20$} as well as selected real-world examples based on the so-called \emph{Penrose square root rule}
and EU population figures. The results allow the estimation of bounds \highlight{for the accuracy of a heuristic solution and can provide informed termination criteria for conventional local search algorithms}. This may be useful in applications where a specific voting power distribution is sought
for a moderate number of council delegates, committee members, or business shareholders. We also analyze the significant magnitude of relative errors that
the mentioned heuristics can produce even for large $n$ in pathological cases.

In the following Section~\ref{sec_voting_systems} we first introduce
binary voting systems and their basic properties. The inverse power index problem is then formalized in Section~\ref{sec_Kurz:2012},
along with a brief discussion of the worst-case bounds which follow from Alon and Edelman's \citeyear{Alon/Edelman:2010} results. In
Section~\ref{sec_methodology} we present the design of our comparative investigation. The corresponding results are the topic of
Section~\ref{sec_results}. We conclude in Section~\ref{sec_conclusion}.

\section{Binary voting systems}
\label{sec_voting_systems}
We consider \emph{binary} voting systems, i.e., each voter $i\in N=\{1,\dots,n\}$ casts a binary vote (e.g., ``yes'' or ``no'')
and this determines a binary collective decision. Such a situation can mathematically be described by a Boolean function
$\chi\colon 2^N\rightarrow \{0,1\}$, where $2^N$ denotes the set of subsets of $N$. A coalition $S\subseteq N$ can, e.g., be
interpreted as the set of ``yes''-voters for a particular (unspecified) proposal.

\begin{definition}\label{def:sg}
  A \emph{simple game} is a monotone Boolean function, i.e., a mapping $\chi\colon 2^N\rightarrow \{0,1\}$ with $\chi(S)\le \chi(T)$
  for all $S\subseteq T$, which additionally satisfies $\chi(\emptyset)=0$, $\chi(N)=1$.
\end{definition}

Coalition $S\subseteq N$ is called \emph{winning} if $\chi(S)=1$, and otherwise losing. $S\subseteq N$ is called a \emph{minimal
winning coalition} if it is winning and all proper subsets are losing. A simple game is uniquely characterized by its set of minimal winning coalitions.
We refer the reader to \citeN{Taylor/Zwicker:1999} \highlight{for a detailed introduction to simple games}.

\highlight{As illustration} consider the simple game $v$ which is characterized by the set $\big\{\{1\},\{2,3\}\big\}$ of its minimal
winning coalitions. Taking all supersets of the minimal winning coalitions yields $\big\{\{1\},\{1,2\},\{1,3\},\{1,2,3\},\{2,3\}\big\}$
as the set of all winning coalitions.\footnote{\highlight{All minimal winning and hence all winning coalitions in a \emph{proper} simple game
have non-empty intersection. So the considered example $v$ is not proper. We do not rule out improper games. They can play a meaningful role even in majoritarian democratic systems, e.g.,
if a sufficiently big parliamentary minority can install a special investigation committee or call for a referendum.}} The remaining subsets are
losing.

The monotonicity imposed in Definition~\ref{def:sg} is a very weak requirement. By introducing Isbell's \emph{desirability
relation}, i.e., $i\sqsupset j$ for two voters $i,j\in N$ if and only if $\chi(\{i\}\cup S\backslash\{j\})\ge \chi(S)$ is satisfied
for all $\{j\}\subseteq S\subseteq N\backslash\{i\}$ (see, e.g., \citeNP{Isbell:1956}), one can define a particularly relevant subclass
of the set of all simple games $\mathcal{S}$:

\begin{definition}
  A simple game $\chi$ is called \emph{complete} (also called \emph{directed}) if the binary relation $\sqsupset$ is a total
  preoder, i.e.,
  \begin{enumerate}
   \item[(1)] $i\sqsupset i$ for all $i\in N$,
   \item[(2)] $i \sqsupset j$ or $j \sqsupset i$ (including ``$i \sqsupset j$ and $j \sqsupset i$'') for all $i,j\in N$, and
   \item[(3)] $i \sqsupset j$, $j \sqsupset h$ implies $i \sqsupset h$ for all $i,j,h\in N$.
  \end{enumerate}
\end{definition}

In our small example one can easily check that $1\sqsupset 2\sqsupset 3$. So $v\in \mathcal{C}$ where
$\mathcal{C}\subset\mathcal{S}$ denotes the set of all complete simple games. Note that also $3\sqsupset2$, i.e., voters~$2$
and $3$ are equally desirable.

Many binary voting systems which are used in practice belong to a further refinement of $\mathcal{S}$:

\begin{definition}
   A simple game $\chi$ is \emph{weighted} if there exist non-negative weights $w_i\in\mathbb{R}_{\ge 0}$ and a positive quota
   $q\in\mathbb{R}_{>0}$ such that $\chi(S)=1$ if and only if $\sum_{i\in S} w_i\ge q$.
\end{definition}

A weighted representation of our small \highlight{illustrative} example is given by $(q;w)=[2;2,1,1]$. We call weighted simple games \emph{weighted voting
games} and denote their collection by $\mathcal{W}$. Every weighted voting game is complete while not every complete simple game
is weighted, i.e., $\mathcal{W}\subset \mathcal{C}\subset\mathcal{S}$.\highlight{\footnote{But each complete simple game and even each simple game can
be represented as the \emph{intersection} of $1\le k<\infty$ weighted voting games. The minimal number $k$ of weighted voting games is
called the \emph{dimension} of the simple game in question (see, e.g., \citeNP{Taylor/Zwicker:1999}, \citeNP{Deineko/Woeginger:2006}).}}
The presently known enumeration results for the three considered classes of binary voting systems are summarized in
Table~\ref{table_number_of_voting_systems} (up to isomorphisms). See, e.g., \citeN{Kurz:2012:representation} for details.

\begin{table}[t]
  \begin{center}
\setlength\extrarowheight{3pt}
   {
    \begin{tabular}{|l|rrrrrrrrrr|}
\hline
      \hspace{0.6cm} $n$&$1$&$2$&$3$& \multicolumn{1}{c}{$4$} & \multicolumn{1}{c}{$5$} &\multicolumn{1}{c}{$6$} &\multicolumn{1}{c}{$7$}
      &\multicolumn{1}{c}{$8$} &\multicolumn{1}{c}{$9$} & \multicolumn{1}{c|}{\highlight{$10$}}\\
      \hline
      $\mathbf{\#}\mathcal{S}$&$1$&$3$&$8$&$28$&$208$&$16351$&$>\!4.7\cdot 10^8$&$>\!1.3\cdot 10^{18}$&$>\!2.7\cdot 10^{36}$ & \highlight{unknown}\\
      $\mathbf{\#}\mathcal{C}$&$1$&$3$&$8$&$25$&$117$&$1171$&$44313$&$16175188$&$284432730174$ & \highlight{unknown}\\
      $\mathbf{\#}\mathcal{W}$&$1$&$3$&$8$&$25$&$117$&$1111$&$29373$&$2730164$&$989913344$ & \highlight{unknown}\\
      \hline
    \end{tabular}}
    \caption{Number of distinct simple games, complete simple games, and weighted voting games}
    \label{table_number_of_voting_systems}
  \end{center}
\end{table}

There are several equivalent representations of binary voting structures besides Boolean functions and lists of minimum winning
coalitions. Simple games can, for instance, be described as independent sets in a graph, and \citeN{Carreras/Freixas:1996} have
introduced a very efficient matrix parameterization of $\mathcal{C}$. Our computation algorithms will exploit yet another possibility.
Namely, we use that voting systems can be represented as points of a polyhedron which have integer coordinates only: for each
$S\subseteq N$ define $0\le x_S\le 1$ 
and add the constraints $x_\emptyset=0$, $x_N=1$, and $x_S\le x_T$ for all $\emptyset\subseteq S\subset T\subseteq N$. Each integer
solution $\highlight{(x_{\emptyset},x_{\{1\}},\ldots, x_N) }\in \{0,1\}^{2^n}$ of this system of linear inequalities
\highlight{corresponds one-to-one} to a simple game (with $\chi(S)=x_S$). Complete
simple games and weighted voting games are described by additional constraints and auxiliary variables \highlight{$w_i\in\mathbb{R}_{\ge 0}$}
for the weights. \highlight{For instance, the inequality
\begin{equation}
  \label{ie_big_M_constraint_weighted1}
  w(S)=\sum_{i\in S} w_i \le q-1 +M\cdot x_S
\end{equation}
implies that the sum of weights of a coalition~$S$ is at most $q-1$
if coalition~$S$ is losing, i.e., $x_S=0$. For a sufficiently large $M$, which can be computed explicitly,
inequality~(\ref{ie_big_M_constraint_weighted1}) is automatically satisfied for $x_S=1$. Similarly the constraint
\begin{equation}
  \label{ie_big_M_constraint_weighted2}
  \sum_{i\in S} w_i \ge q+M(1-x_S)
\end{equation}
forces winning coalitions to have a weight sum exceeding or meeting the quota.}
\highlight{Because each weighted voting
game admits a representation $(q; w_1, \ldots, w_n)$ where $w(S)\le w(T)-1$  for all losing coalitions~$S$ and all winning coalitions~$T$, inequalities (\ref{ie_big_M_constraint_weighted1}) and (\ref{ie_big_M_constraint_weighted2}) capture the weightedness requirement for a simple game described by $(x_{\emptyset},x_{\{1\}},\ldots, x_N)$.} 

\section{The inverse power index problem}
\label{sec_Kurz:2012}
Power indices are mappings from a set of feasible voting structures, such as $\mathcal{S}$ or $\mathcal{W}$, to non-negative real
vectors which are meant to quantify the influence of the members of a voting body on collective decisions. The inverse power index
problem consists in finding a voting system, e.g., $(q;w)\in \mathcal{W}$, which induces a power distribution as close as possible
to a desired one. More formally, for a given number $n$ of voters, the general \emph{inverse power index problem} involves a set
$\Gamma$ of feasible voting structures for $n$ players, a power index $\phi\colon \Gamma\rightarrow\mathbb{R}^n_{\ge 0}$, a desired
power distribution $\beta\in\mathbb{R}^n_{\ge 0}$, and a metric $d\colon \mathbb{R}^n\times
\mathbb{R}^n\rightarrow \mathbb{R}_{\ge 0}$ which measures the deviation between two power vectors. Of course, $d(x,y)=\Vert x-y\Vert$
is a suitable choice for any vector norm $\Vert\cdot\Vert$. Given these ingredients the inverse power index problem amounts to
finding a solution to the minimization problem
\begin{equation}\label{eq:min_problem_stated}
\min_{\chi\in\Gamma}\,d\big(\phi(\chi),\beta\big).
\end{equation}

In this paper, we consider the special instances of this problem where $\Gamma\in \{\mathcal{S},\mathcal{C},\mathcal{W}\}$. We
include $\mathcal{S}$ and $\mathcal{C}$  because they are significantly larger domains for $n\ge 5$ (see
Table~\ref{table_number_of_voting_systems}) and some prominent real-world electoral systems fail to correspond to weighted voting games.
Examples include the current voting rules (Treaty of Nice) and the future ones (Treaty of Lisbon) for the EU Council, which \highlight{require
majorities in more than one dimension (e.g., the Nice rules call for 255 out of 345 votes, 14 out of 27 member states, and 62\% of EU population).}
We take the (normalized) Penrose-Banzhaf index $B(v)$ as the voting power index of interest.

\begin{definition}
For a given $n$-player simple game $\chi$ the \emph{absolute Penrose-Banzhaf index} $B'_i(\chi)$ for player $i$ is defined as
$$
  B'_i(\chi)=\frac{1}{2^{n-1}}\cdot\sum_{\emptyset\subseteq S\subseteq N\backslash\{i\}} \chi(S\cup\{i\})-\chi(S).
$$
The \emph{(normalized) Penrose-Banzhaf index} (PBI) $B_i(\chi)$ for player~$i$ is defined as
$$
B_i(\chi)=\frac{B'_i(\chi)}{\sum_{j=1}^nB'_j(\chi)}.
$$
\end{definition}

Our distance computations will be based on the $\Vert \cdot \Vert_1$-norm (i.e., the sum of deviations between $B_i(\chi)$ and
$\beta_i$ for all players~$i$), the $\Vert \cdot \Vert_\infty$-norm (i.e., the maximum deviation), and a weighted version of the
former. Section~\ref{sec_methodology} will provide more details.

To the best of our knowledge, there exists only one (non-trivial) non-approxi\-ma\-tive result on how well the inverse problem
can be solved for the PBI in the worst case. \highlight{For completeness and later reference we include this rather recent finding by
\citeN{Alon/Edelman:2010} here. It considers a given game $v$ with $n$ players in which $1-\eps$ of the total (normalized) PBI
is concentrated amongst $k<n$ ``major'' players. Alon and Edelman then provide a construction for a game
$\tilde{v}$ such that the worth $\tilde{v}(S)$ of a coalition~$S$ depends only on $T=S\cap\{1,\dots,k\}$, i.e., the subset of major
players in $S$ fully determines whether $S$ is winning in game $\tilde v$ or not.\footnote{\highlight{Specifically,
one sets $\tilde{v}(S)=1$ if and only if $\sum_{U\subseteq\{k+1,\dots,n\}} v(T\cup U)\ge \frac{2^{n-k}}{2}$. So the coalition
$T$ of major players~-- and hence all supersets $S'=T\cup U$ that are obtained by adding different coalitions $U$ of ``minor'' players~-- is
winning in $\tilde v$ if and only if a majority of the latter coalitions $S'$ are winning in $v$.}} One can easily observe that $\tilde{v}$
is a simple game
and $B_i(\tilde{v})=0$ for all $i>k$.\footnote{\highlight{Moreover, we have $\tilde{v}\in\mathcal{C}$ if
$v\in\mathcal{C}$ and $\tilde{v}\in\mathcal{W}$ if $v\in\mathcal{W}$, i.e., the construction respects completeness or weightedness of the
simple game in question.}} The essential finding of Alon and Edelman then is that the deviation $\Vert B(v)-B(\tilde{v})\Vert_1$
is bounded from above by a function which depends on $\varepsilon$ and the number $k$ of major players only. Considering the $k$-player
simple game $v'$, which arises from $\tilde{v}$ by removing the $n-k$ null players, their result can be stated as follows:}

\begin{theorem}[Alon-Edelman]
  \label{thm_alon}
\highlight{Consider the simple game $v$ with players $N=\{1, \ldots, k,$ $\ldots, n\}$ and  $0<\eps<\frac{1}{k+1}$.} If $\sum_{i=k+1}^nB_i(\chi)\le\varepsilon$,
then there exists a simple game $\chi'$ with $k$ voters such that
$$
    \sum_{i=1}^k \left|B_i(\chi)-B_i(\chi')\right|\,+\,\sum_{i=k+1}^n B_i(\chi)
    \le\frac{(2k+1)\varepsilon}{1-(k+1)\varepsilon}+\varepsilon.
$$
\end{theorem}

\highlight{This result is very useful for obtaining lower bounds on distances in the context of the inverse problem because one may suitably reduce the problem from $n$
to $k$~players and only make an error with the indicated bound.
Specifically, let $\beta=(\beta_1,\dots,\beta_n)$ be a desired power distribution and $k,\varepsilon$ be constants satisfying $0<\eps<\frac{1}{k+1}$. Let us denote the unknown $n$-player simple game
whose PBI has smallest distance to $\beta$ by $v^*$. One can then bound $\Vert B(v^*)-\beta \Vert_1$ by distinguishing two cases.}

\highlight{First, suppose that $\sum_{i=k+1}^n B_i(v^*)\ge \varepsilon$. In this case, we can only apply Theorem~\ref{thm_alon} if the inequality happens to be tight. However, we know that the deviations between $B_i(v^*)$ and $\beta_i$ for players $k+1, \ldots, n$ are at least as big as
$$
\big| \sum_{i=k+1}^n B_i(v^*) - \sum_{i=k+1}^n \beta_i \big|.
$$
This can range from $1 - \sum_{i=k+1}^n \beta_i$ when $\sum_{i=k+1}^n B_i(v^*)=1$ to $|\eps - \sum_{i=k+1}^n \beta_i|$ when $\sum_{i=k+1}^n B_i(v^*)$ $= \varepsilon$ in the considered case. Analogous reasoning applies to the additional deviations for players $1, \ldots, k$, and we thus have the bound
$$
  \Vert B(v^*)-\beta\Vert_1\ge \min_{1\ge x\ge\varepsilon} \bigl|1-x-\sum_{i=1}^k \beta_i\bigr|+\bigl|x-\sum_{i=k+1}^n \beta_i\bigr|=:l_1.
$$
}

\highlight{In the second case, i.e., when $\sum_{i=k+1}^n B_i(v^*)< \varepsilon$, Theorem~\ref{thm_alon} applies. It tells us that there exists some $k$-player simple game $v'$ with
$$
\sum_{i=1}^k \left|B_i(v^*)-B_i(\chi')\right|\,+
\,\sum_{i=k+1}^n \left|B_i(v^*)-0\right| \le\frac{(2k+1)\varepsilon}{1-(k+1)\varepsilon}+\varepsilon.
$$
If we use $\tilde v$ to denote the $n$-player game that extends $v'$ by adding $k+1, \ldots, n$ as null players (so that $B(\tilde v)=(B_1(v'), \ldots, B_k(v'), 0, \ldots, 0)$), this can also be written as
\begin{equation}\label{eq:AEfortildev}
\Vert B(v^*)-B(\tilde v)\Vert_1  \le\frac{(2k+1)\varepsilon}{1-(k+1)\varepsilon}+\varepsilon.
\end{equation}
Now if we solve the $k$-player inverse problem for the (typically non-normalized) $k$-vector $\beta'=(\beta_1,\dots,\beta_k)$ which coincides with the first $k$ components of $\beta$, then the resulting minimal distance $\varepsilon'$ is a lower bound for the distance between $B(\tilde v)$ and $\beta$, i.e.,
\begin{equation}\label{eq:epsprime}
\Vert B(\tilde v)-\beta\Vert_1\ge \varepsilon'.
\end{equation}
}

\highlight{We can then appeal to the triangle inequality for metric $d_1(x,y)=\Vert x-y \Vert_1$ and conclude
\begin{eqnarray*}
\Vert B(v^*)-\beta \Vert_1 \ge \Vert B(\tilde v)-\beta \Vert_1  - \Vert B(\tilde v)-B(v^*)\Vert_1 \ge \varepsilon' -  \frac{(2k+1)\varepsilon}{1-(k+1)\varepsilon}+\varepsilon:=l_2
\end{eqnarray*}
from (\ref{eq:AEfortildev}) and (\ref{eq:epsprime}).
Thus, in either of the two cases we have $\Vert B(v^*)-\beta\Vert_1\ge \min(l_1,l_2)$.}

\medskip

\highlight{Let us illustrate this by an example and} suppose that one seeks to find a voting game with a power distribution as close as possible to
$\beta^n=(0.75,0.25,0,\ldots,0)\in \mathbb{R}^n_{\ge 0}$ for $n\ge 2$.
\highlight{We will show $\Vert B(v)-\beta^n\Vert_1\ge\frac{1}{9}$ for \emph{all} $n$-player simple games $\chi$.%
\footnote{\highlight{Such
artificially constructed target vectors $\beta^n$, for which Alon and Edelman's results have bite, may not be of much practical
relevance. But they indicate the problems of
requiring any fixed level of accuracy in the stopping rule of a local search algorithm. Moreover, one can conceive of real-world enterprises in which the majority and minority partners indeed seek to split
voting power $3:1$ and want to render all $n-2$ other stakeholders null players.}}
To this end we choose $k=2$
and $\varepsilon=\frac{1}{18}$. If $\sum_{i=3}^n B_i(v^*)\ge\varepsilon$ for the unknown distance-minimizing game $v^*$ then we have
$$
  \min_{1\ge x\ge\varepsilon} \bigl|1-x-\sum_{i=1}^k \beta_i\bigr|+\bigl|x-\sum_{i=k+1}^n \beta_i\bigr|
  =\min_{1\ge x\ge\varepsilon} 2x=2\varepsilon=\frac{1}{9}.
$$
So in this case we have $\Vert B(v^*)-\beta^n\Vert_1\ge\frac{1}{9}$. In the other case of $\sum_{i=3}^n B_i(v^*)<\varepsilon$, we solve the inverse power index problem for $k=2$ players and
$\beta'=\beta^2=(0.75,0.25)$. Since the only possible 2-player PBI vectors are given by $\big\{(1,0), (\frac{1}{2},
\frac{1}{2}), (0,1)\big\}$ we have a minimal deviation of $\varepsilon'=\frac{1}{2}$. Because $\frac{(2k+1)\varepsilon}
{1-(k+1)\varepsilon}+\varepsilon=\frac{7}{18}$ we have $l_2=\frac{1}{9}$ and again conclude $\Vert B(v^*)-\beta^n\Vert_1\ge\frac{1}{9}$.}
Hence, $\beta^n$ cannot be approximated by the PBI of a
simple game with an $\Vert\cdot\Vert_1$-error less than $\frac{1}{9}$. The latter is the sharpest possible bound obtainable from
Theorem~\ref{thm_alon}. It can be improved computationally to slightly more than $\frac{14}{37}$ for $n\le 11$ on $\mathcal{S}$
and for $n\le 16$ on $\mathcal{C}$ and $\mathcal{W}$ (see \citeNP{Kurz:2012}).

\section{Design of the computational investigation}
\label{sec_methodology}
When the inverse problem arises in political applications of constitutional design, PBI vectors $\beta$ which are proportional to
the square root of a population size vector $p$ play an elevated role. The reason is that~-- under the probabilistic assumptions
which underlie the PBI~-- a binary voting system $v$ with $B(v) = \beta$ and
\begin{equation}\label{eq:SQRbeta}
  \beta_i=\frac{\sqrt{p_i}}{\sum_{j=1}^n\sqrt{p_j}}
\end{equation}
would equalize the voting power of citizens in a \emph{two-tier system} in which $n$ delegates adopt the bottom-tier majority
opinion of their respective constituency~$i\in\{1, \ldots, n\}$ and then cast a $w_i$-weighted vote in a top-tier assembly (e.g.,
the EU Council). See \citeN{Penrose:1946}, \citeN{Felsenthal/Machover:1998}, \citeN{Kaniovski:2008} or \citeN{Kurz/Maaser/Napel:2012}
for details. In our computations we will consider this \emph{Penrose square root rule} for varying $n$ and some historical population figures in order
to select target vector \highlight{examples which have a specific political motivation}.

In principle, however, \emph{any} vector in $\mathbb{R}^n_{\ge 0}$ whose entries sum up to $1$ might be a desired power distribution
$\beta$. For instance, partners of a non-profit R\&D joint venture might have made relative financial contributions of
$\big(\frac{1}{3}, \frac{1}{3}, \frac{1}{9}, \frac{1}{9}, \frac{1}{9}\big)$ and possibly want to align a~priori voting power
\highlight{in the directorate} to this vector as well as possible. Ideally, for a given number $n$ of voters, one would like to compare
the exact and heuristic solutions to the inverse problem for \emph{all} possible normalized target vectors $\beta \in \Delta(n-1)$, where
$\Delta(n-1)$ denotes the $n-1$-dimensional unit simplex. \highlight{This is computationally infeasible. We, however, complement our analysis
of politically motivated square root vectors by vectors $\beta$ from a discrete subset of $\Delta(n-1)$, namely a finite grid on $\Delta(n-1)$
with step size 0.01. We also resort to approximations of the exact solution when $n$ is too large.}

We will compare the (approximated) exact solution of the inverse problem on domain $\mathcal{S}$, $\mathcal{C}$, or $\mathcal{W}$
for a given desired PBI $\beta$ with three different heuristics. \highlight{These} stay in the class
$\mathcal{W}$ of weighted voting games and have in common that voting weights are set equal to the desired voting power, i.e., $w=\beta$.
They pick a distinct quota, and hence typically a different voting system $v\in\mathcal{W}$.

The first heuristic~-- referred to as the \emph{50\%-heuristic}~-- just chooses $q^\circ=\frac{1}{2}$. Simple majority is arguably
the most common majority rule in practice. The 50\%-heuristic simply picks it and ignores the potentially large discrepancies between
voting weight and voting power that can arise. \highlight{This can be motivated by the Penrose limit theorem when $n$ is at least moderately big (see fn.~\ref{fn:oceanicgames} however).}

The second, more sophisticated heuristic has been suggested by S\l{}omczy\'nski and \.Zyczkowski
\citeyear{Slomczynski/Zyczkowski:2006,Slomczynski/Zyczkowski:2007}. Their motivation was to implement PBI vectors proportional to the
square root of population sizes in the European Union, but the heuristic applies to arbitrary target vectors. Namely,
the \emph{$q^*$-heuristic} selects the quota
$$
q^*=\frac{1}{2}\cdot\Big(1+\sqrt{\sum_i w_i^2}\Big)
$$
for an arbitrary $w=\beta\in \Delta (n-1)$. \citeN{Slomczynski/Zyczkowski:2007} derive this quota by considering the random weight $W$ which
is accumulated if all coalitions $S\subseteq N$ are \emph{equiprobable}, as the PBI's probabilistic justifications suppose.
Equiprobability at the level of coalitions is equivalent to assuming that each voter~$i\in\{1, \ldots, n\}$ joins the formed coalition
independently of the others with probability $\frac{1}{2}$. The mean of $W$ hence is $\mu=\sum_{i=1}^n\frac{1}{2}w_i=\frac{1}{2}$ and
its variance is $\sigma^2=\frac{1}{4}\sum_{i=1}^n w_i^2$. Being the sum of independent bounded random variables, $W$ is approximately
normally distributed \emph{if} $n$ is sufficiently large and each of the weights is sufficiently small.\footnote{A key technical requirement
is that \highlight{$w_j \ll \sqrt{\sum w_i^2}$ for all $j\in N$, i.e., $w_j\sqrt{n}$ is sufficiently smaller than the quadratic mean of the
weights $\sqrt{\frac{1}{n}\sum w_i^2}$. \label{fn:SZcondition}}} Assuming that this is the case and, therefore, that the discrete random variable
$W$ can be replaced by the continuous one $\tilde W$, the inflection point of the corresponding normal density of $\tilde W$ is located at
$q^*=\mu+\sigma$. Since the second derivative of $\tilde W$'s density vanishes at $q^*$, one can approximate the density by a linear function
with reasonably high accuracy. This linear approximation then allows to establish approximate proportionality of $B(q^*;w)$ and $w$. We refer to \citeN{Slomczynski/Zyczkowski:2007} for details.

Our final heuristic, which we will refer to as the \emph{$\bar q$-heuristic}, replaces $q^*$ by
$$
\bar q=\frac{1}{2}+\frac{1}{\sqrt{\pi n}}.
$$
This quota approximates the expected value of $q^*$ when $\beta$ is proportional to the component-wise square root of a
population size vector $p=(p_1, \ldots, p_n)$ which is drawn from a flat Dirichlet distribution (see \citeNP{Slomczynski/Zyczkowski:2011}).
The motivation for computing such an average is the following: even though the $q^*$-heuristic can approximate the Penrose square root
rule (\ref{eq:SQRbeta}) very transparently for a given population distribution $p$, frequent changes in the population would call not only
for frequent changes of the \highlight{prescribed} voting weights $w$ but also of the quota $q^*$. That current voting weights in the
EU already recur to population figures, which are updated on an annual basis, suggests that weight changes may be regarded as unproblematic. A varying
decision threshold~-- perhaps $q=65\%$ in one year, $q'=61\%$ in the next, then $q''=67\%$, etc.~-- however seems politically less palatable. It
may then make sense to average $q^*$ over a wide range of values for $w=\beta\propto \sqrt{p}$, and the $\bar q$-heuristic simply assumes
that all population distributions $p\in \Delta(n-1)$ are equally likely.\footnote{The expected value of the $p$-specific optimal quotas
$q^*(p)$ for a particular (e.g., Dirichlet) distribution of $p$, of course, need not coincide with the quota that is optimal when $p$ is
treated as a random variable. Stochastic optimization techniques are likely to yield a somewhat better $q$-heuristic than the one suggested
by \citeN{Slomczynski/Zyczkowski:2011}.} Because $\bar q \to \frac{1}{2}$ as $n\to \infty$, the 50\%-heuristic is the limit of the
$\bar q$-heuristic and can be viewed as an approximation of it for not too small $n$.

Let us remark that investigations by \citeN{Kurth:2008} have called attention to numerical problems when heuristics which involve
irrational voting weights and quotas, as the $q^*$ or $\bar q$-heuristics commonly do, are implemented.
Rounding after, e.g., 4 decimal places 
can result in voting systems which differ significantly from what was intended. Because it is impractical to deal with weights of a
hundred decimal places or more, it is attractive to work with the underlying Boolean functions or integer points of a suitable polyhedron
as long as possible, and to determine minimal integer weights $w$ and a quota $q$ which \highlight{efficiently} represent a given
$v\in \mathcal{W}$ when needed.\footnote{A \emph{minimal} integer representation of a weighted voting game has the 
advantage that the PBI and other power indices can be computed 
\highlight{particularly quickly}.} We use this approach
here whenever possible, and refer the interested reader to \highlight{Freixas and Molinero  \citeyear{Freixas/Molinero:2009,Freixas/Molinero:2010}},
\highlight{\citeN{Freixas/Kurz:2011}}, or \citeN{Kurz:2012:representation}.

We calculate the globally optimal solution to the inverse problem for a given target PBI $\beta$ by complete enumeration of the elements in
the respective class of binary voting systems for $n\le 7$ (see Table~\ref{table_number_of_voting_systems}). For larger $n$, we mostly
focus on approximations of the exact solution. \highlight{These are obtained either by a hill-climbing algorithm or, preferably, by ILP techniques.
How the latter are used is explained in the Appendix in detail. The implemented ILP-based bisection algorithm would yield globally optimal solutions
when given enough running time and memory. We interrupted it for efficiency reasons whenever a desired precision had been reached.} The key
idea \highlight{of the ILP-based approach} is to consider the integer polyhedron which contains all simple games whose PBI is less than a
given factor $\alpha>0$ away from the desired vector $\beta$. \highlight{It can be checked by using readily available ILP solver software if}
this polyhedron is empty. \highlight{Then,} no such game exists and $\alpha$ needs to be raised. If not, $\alpha$ can be lowered. The minimal
level of $\alpha$ (or an approximation with desired precision) together with the corresponding voting systems, can thus be found by the bisection
method: \highlight{namely, by iteratively halving the interval defined by the best lower and upper bounds that have been computed so far. Pseudo-code of the algorithm is provided in Appendix~A.}

In evaluating the quality of the mentioned heuristics, we consider distances to the desired power vector, $\beta$, and to the globally
optimal one, $B(v^*)$, in three different metrics. The first one is the metric $d_1(x,y)=\Vert x-y\Vert_1=\sum_{i=1}^n |x_i-y_i|$ induced
by the $\Vert\cdot\Vert_1$-norm, which is also considered in Theorem~\ref{thm_alon}. The second is the metric induced by the
$\Vert\cdot\Vert_\infty$-norm, i.e., $d_\infty(x,y)=\Vert x-y\Vert_\infty=\max_{i\in\{1,\ldots, n\}} |x_i-y_i|$. We refrain from also
considering the Euclidean metric induced by the $\Vert\cdot\Vert_2$-norm, which has been considered, e.g., by \citeN{Slomczynski/Zyczkowski:2007}.
The reason is that this would turn the ILP formulation of the inverse power index problem into a \emph{binary non-linear programming} one.
This would be considerably harder to solve and add relatively little information because $\Vert x\Vert_\infty\le\Vert x\Vert_2\le\sqrt{n}
\Vert x\Vert_\infty$ for all $x\in\mathbb{R}^n$.

More interesting, in our view, is a variation of $d_1$ which takes the Bernoulli model that underlies the PBI and Penrose's square root
rule seriously. This model assumes that all bottom-tier voters in constituency~$i\in \{1, \ldots,n\}$ cast a ``yes'' or ``no'' vote
equiprobably and independently of all others. The probability for one of $p_i$ individual voters in constituency~$i$ to be \emph{pivotal}
for the constituency's aggregate decision~-- i.e., to induce the $i$-delegate at the top-tier council to cast voting weight $w_i$ in favor
of ``yes''  by individually voting ``yes'', and ``no'' by voting ``no''~-- is approximately $\sqrt{2/(\pi p_i})$. The joint probability of
a given voter being pivotal in his \highlight{or her} constituency~$i$ \emph{and} of this constituency being pivotal at the top tier is hence $B_i(\chi)\cdot
\sqrt{2/(\pi p_i})$. This is why the square root PBI vector in equation (\ref{eq:SQRbeta}) equalizes the indirect influence of citizens on
collective decisions across constituencies. If one now weights any deviation between \highlight{(i)} the probability for a given voter in constituency~$i$
to be doubly pivotal and \highlight{(ii)} the egalitarian ideal of $\beta_i\cdot \sqrt{2/(\pi p_i})$ with $\beta_i=\sqrt{p_i}/\sum_{j=1}^n \sqrt{p_j}$
equally, then the total misrepresentation associated with the top-tier voting system $v$ amounts to $$
\sum_{i=1}^n p_i \cdot \big|\beta_i - B_i(v) \big|\cdot \sqrt{2/(\pi p_i})=c\cdot \sum_{i=1}^n \sqrt{p_i} \cdot \big|\beta_i - B_i(v) \big|
$$
for $c>0$. Whenever the desired vector $\beta$ is derived from Penrose's square root rule and a vector $p$ which represents EU population data,
we will, therefore, also consider the variation of metric $d_1$ which weights absolute deviations by the square root of relative
population, i.e., study the metric\footnote{Consideration of a similar variation of $d_\infty$ broadly confirms the comparisons based on
$d_1$, $d_{1}'$, and $d_\infty$.}
$$
d_{1}'(x,y)= \sum_{i=1}^n \sqrt{\frac{p_i}{\sum_{j=1}^n p_j}}\cdot|x_i-y_i|.$$

\section{Computational results}
\label{sec_results}

In this section we present our numerical results. Subsection~\ref{subsec:EU} considers the EU Council of Ministers as a prototype
of a real-world weighted voting system. We then look at the entire \highlight{discretized} space of possible power distributions for
$n\le 7$ and  random samples thereof in Subsection~\ref{subsec:grid}. In order to study analytically how deviations between simple heuristics
and actual optimization depend on $n$, we investigate a particular parametric example in Subsection~\ref{subsec:analytical}.

\highlight{Exact solutions to the inverse power index problem
that are reported in Tables~\ref{table_quality_d1}--\ref{table_quality_d2} have been obtained using our ILP-based bisection algorithm,
as described above and in Appendix~A. For Table~\ref{table_grid_optimal} we have used exhaustive enumeration of all possible Banzhaf vectors
for $n\le 7$,  and resorted to approximations obtained by hill-climbing algorithms for $8\le n\le 20$. The bisection part of the ILP approach was implemented in C++ while we used
the ILOG CPLEX Interactive Optimizer 12.4.0.0 in order to solve individual integer linear programming problems. The hill-climbing algorithms used in
Subsection~\ref{subsec:grid} were implemented in C++. We employed a Quad-Core AMD Opteron processor with 2700~Mhz, 132~GB RAM, and
a cache size of 512~KB on a 64-bit Linux system as our hardware.}

\subsection{Examples of real-world weighted voting systems}\label{subsec:EU}

We first consider the (EEC or EC or)
EU Council of Ministers in the years 1958, 1973, 1981, 1986, 1995, 2006, and 2011 with respectively $n\in\{6,9,10,12,15,25,27\}$
members as examples. The historical population data for $n\in\{6,\ldots,15\}$ are taken from \citeN[sec.~5.3]{Felsenthal/Machover:1998}, the
data for $n\in\{25,27\}$ are official Eurostat figures downloaded on 19.01.2012. \highlight{The target power distributions $\beta$ are the respective ``fair'' ones}
computed by Penrose's square root rule (see equation~(\ref{eq:SQRbeta})).

In Tables~\ref{table_quality_d1}--\ref{table_quality_d2} we compare the three considered heuristics under different metrics
with the \highlight{(approximated)} optimal solution
of the inverse power index problem. \highlight{We distinguish between $\mathcal{S}$, $\mathcal{C}$, and $\mathcal{W}$ as the sets of admissible
voting structures. Besides the absolute deviations (measured in the respective metric) we also report an indicator of relative quality: if the
distance between $\beta$ and the PBI $B(v^*)$ of the optimal solution $v^*\in \mathcal{S}$ is $\alpha$, then
this is the \emph{unavoidable} absolute ``error'' associated with the given instance of the inverse problem. Now if a certain heuristic
delivers a distance of $\delta$ then $({\delta-\alpha})/{\alpha}$ can be regarded as the \emph{avoidable error}
relative to global optimization in $\mathcal{S}$. It is labeled \emph{$\mathcal{S}$-error} in the
tables.} A value of 1 \highlight{(or 100\%)} means that the heuristic's approximation error is twice the unavoidable one.

The ``$\dagger$''-symbol indicates that the stated value in Tables~\ref{table_quality_d1}--\ref{table_quality_d2} has not been
computationally proven to be optimal: for simple games and $n=9$, for instance, we stopped the ILP solution process after memory
usage of 31~GB and 18461700 branch-and-bound nodes; for $n=10$, 
we interrupted after 301~GB and 16735508 nodes. 
\highlight{The ``$\dagger\dagger$''-symbol indicates that a lower bound for
the minimal distance in $\mathcal{S}$ or $\mathcal{C}$ was inferred from $\mathcal{W}$. The ``$\dagger\dagger$''-marked numbers need
not be optimal a fortiori. The ``$\infty$'' entries indicate avoidable errors greater than factor $500$. Finally, $0.000000^\dagger$
or $0.000000^{\dagger\dagger}$ represent positive numbers $<0.5\cdot 10^{-6}$.}

\highlight{The computation times for obtaining the
numbers in the $\mathcal{C}$-column in Table~\ref{table_quality_d1}, using the hardware and software described above, ranged from
less than a second for $n=6$ to 5~days for $n=12$. The exact solution for $n=10$ in $\mathcal{W}$ took 2~days; the approximate
one for $n=27$ was obtained in 3~hours.}\footnote{\highlight{The additional constraints (\ref{ie_big_M_constraint_weighted1})
and (\ref{ie_big_M_constraint_weighted2}) which ensure $v$'s weightedness in our ILP formulation considerably slow down the
computations (because the so-called \emph{integrality gap} increases). In contrast, the analogous constraints which ensure
completeness impose useful structure on the problem compared to unrestrained optimization in $\mathcal{S}$. This explains, e.g., why
an exact solution can be reported for $n=10$ in Table~\ref{table_quality_d2} in the $\mathcal{C}$-column but not in the $\mathcal{S}$
and $\mathcal{W}$-columns.}}

\begin{table}[htp]
  \begin{center}
    \begin{tabular}{|r|l|l|l|r@{}r|r@{}r|r@{}r|}
      \hline
       & \multicolumn{1}{c|}{$v^*\in \mathcal{S}$}  & \multicolumn{1}{c|}{$v^{**}\in \mathcal{C}$} & \multicolumn{1}{c|}{$v^{***}\in \mathcal{W}$}
       &\multicolumn{2}{c|}{50\%-heuristic} & \multicolumn{2}{c|}{$q^*$-heuristic}
      & \multicolumn{2}{c|}{$\bar{q}$-heuristic} \\
      \multicolumn{1}{|c|}{$n$} & \multicolumn{1}{c|}{$d_1$} & \multicolumn{1}{c|}{$d_1$} & \multicolumn{1}{c|}{$d_1$} & \multicolumn{1}{c}{$d_1$}
      & $\mathcal{S}$-error & \multicolumn{1}{c}{$d_1$} & $\mathcal{S}$-error & \multicolumn{1}{c}{$d_1$} & $\mathcal{S}$-error \\
      \hline
    6 &0.051857 & 0.051857 & 0.051857 & 0.300398 & 4.79    & 0.091100 & 0.76  & 0.091100 & 0.76   \\
    9 &0.005294$^\dagger$  & 0.008641 & 0.010359  & 0.065528 & 11.38   & 0.060195 & 10.37 & 0.069792 & 12.18  \\
   10 &0.002639$^\dagger$& 0.004840   & \highlight{0.007219} & 0.038751 & 13.68   & 0.033229 & 11.59 & 0.026466 & 9.03   \\
   12 &\highlight{0.001033$^\dagger$}& \highlight{0.001033$^\dagger$}& 0.005170$^\dagger$& 0.028700 & \highlight{26.78}   & 0.019827 &
   \highlight{18.19}  & 0.019827 & \highlight{18.19}\\
   15 & 0.000476$^{\dagger\dagger}$ & 0.000476$^{\dagger\dagger}$ & 0.000476$^\dagger$& 0.026742 & 55.18   & 0.006820 & 13.33 & 0.006361 & 12.36  \\
   25 & 0.000000$^{\dagger\dagger}$ & 0.000000$^{\dagger\dagger}$ & 0.000000$^\dagger$& 0.019422 & ``$\infty$''      & 0.000744 & ``$\infty$''
   & 0.003096 & ``$\infty$''     \\
   27 & 0.000000$^{\dagger\dagger}$ & 0.000000$^{\dagger\dagger}$ & 0.000000$^\dagger$& 0.018003 & ``$\infty$''      & 0.000633 & ``$\infty$''
   & 0.002457 & ``$\infty$''     \\
      \hline
    \end{tabular}
    \caption{Performance for Penrose square root targets in the $d_1$-metric (1958--2011 EU data)}
    \label{table_quality_d1}
  \end{center}
\end{table}

\highlight{The reported numbers give rise to several observations that are independent of the chosen metric:\footnote{Note that the three metrics behave differently when, e.g., distance between
$(1,0,\ldots,0)$ and $(\frac{1}{n}, \ldots, \frac{1}{n})\in \Delta(n-1)$ is considered for increasing $n$. Deviations
should, therefore, be compared only within and not across tables.}} (i) The approximation errors of
the heuristics and the optimal solutions in $\mathcal{W}$ (and a~fortiori in $\mathcal{C}$ and $\mathcal{S}$) tend to zero
as $n$ increases. (ii) Except for $n=9$, the $q^*$-
and the $\bar{q}$-heuristics perform noticeably better than the 50\%-heuristic. \highlight{(iii) The $q^*$ and $\bar q$-heuristics
produce comparable absolute deviations from the ideal for $n\le 15$ but differ by a factor of 2 or more for $n>15$. We conjecture that this has to do with the
normal density approximation, which is underlying \citeANP{Slomczynski/Zyczkowski:2007}'s \citeyear{Slomczynski/Zyczkowski:2007}
derivation of $q^*$, becoming noticeably more accurate for the population distribution in the enlarged EU. This allows
$q^*$'s performance to improve by almost an order of magnitude between $n=15$ and $25$, while performance of the $\bar q$-heuristic
(which ignores the specific population distribution at hand and picks the optimal quota averaged over many possible distributions)
improves only by about the same factor as the $50\%$-heuristic and thus in line with the asymptotic proportionality results captured by the
Penrose limit theorem.}

\highlight{We also find that (iv)} the respective optimal weighted games $v^{***}\in \mathcal{W}$ yield deviations that are only moderately
higher than those of $v^*\in \mathcal{S}$ \highlight{in absolute terms. This might be interpreted as indicating that relatively little is lost by restricting attention to weighted voting
games in conventional hill-climbing algorithms. It should be noted, however, that we could not prove optimality in $\mathcal{S}$ for $n\ge 9$, and for $n\ge 15$ we only have upper bounds obtained from $\mathcal{W}$. So the observation might not be very robust. And, in relative terms, the errors in $\mathcal{S}$ or $\mathcal{C}$ are several times smaller than those in $\mathcal{W}$ for $n=10$ or 12.}

\highlight{Finally, (v) the \emph{relative} errors of the heuristics compared to either
$v^*\in \mathcal{S}$ or $v^{***}\in \mathcal{W}$ are sizeable even for small $n\le 15$; and they become huge for $n>25$. For small $n$ like
$n=6$ or 9, the unavoidable error, i.e., the distance between $B(v^*)$ and the target vector $\beta$, is still big because comparatively
few distinct PBI values exist. At the same time, such numbers $n$ are far too small for the normal approximation which underlies the
$q^*$ and $\bar q$-heuristics or for the asymptotics which motivate the 50\%-heuristic to have leverage. So the heuristics do not perform
well in absolute terms, but they are not that bad in relative terms because of high unavoidable errors. Now as $n$ increases, the
heuristics perform significantly better in absolute terms. However, the unavoidable error vanishes even more quickly as the number
of distinct simple voting games and, hence, of feasible PBI vectors increases very fast in $n$ (see Table~\ref{table_number_of_voting_systems}).
}

Observation (v) is probably the most interesting: whenever one seeks an optimal solution of the inverse
power index problem, all three heuristics are unsatisfactory from a pure operations research perspective. The heuristic
solutions can be improved by very large factors, and this becomes more rather than less pronounced as $n$ grows. Of course, from an applied point of view the absolute approximation errors get so small for large $n$ that
they may be regarded as negligible. They might still be relevant, however. \highlight{To get a sense} for what a deviation at the 5th decimal place means consider, e.g., the ideal Penrose square root power
distribution $\beta^{27}$ for the EU Council from 2011 and \highlight{compute} the analogous vector $\beta^{27\prime}$ which would result if
50000 people moved from Germany to France or were mis-counted in the statistics. Then  
$\Vert\beta^{27}-\beta^{27\prime}\Vert_1\approx 0.0000634$. 

\begin{table}[t]
  \begin{center}
    \begin{tabular}{|r|l|l|l|r@{}r|r@{}r|r@{}r|}
      \hline
       & \multicolumn{1}{c|}{$v^*\in \mathcal{S}$}  & \multicolumn{1}{c|}{$v^{**}\in \mathcal{C}$} & \multicolumn{1}{c|}{$v^{***}\in \mathcal{W}$}
       &\multicolumn{2}{c|}{50\%-heuristic} & \multicolumn{2}{c|}{$q^*$-heuristic}
      & \multicolumn{2}{c|}{$\bar{q}$-heuristic} \\
      \multicolumn{1}{|c|}{$n$} & \multicolumn{1}{c|}{$d_1'$} & \multicolumn{1}{c|}{$d_1'$} & \multicolumn{1}{c|}{$d_1'$} & \multicolumn{1}{c}{$d_1'$}
      & $\mathcal{S}$-error & \multicolumn{1}{c}{$d_1'$} & $\mathcal{S}$-error & \multicolumn{1}{c}{$d_1'$} & $\mathcal{S}$-error \\
      \hline
       6 & 0.018967\hspace{4pt}        & 0.021487\hspace{4pt}        & 0.021487\hspace{4pt}        & 0.110284 & 4.81  & 0.027465 & 0.45 & 0.027465 & 0.45 \\
       9 & 0.001902$^\dagger$& 0.002752\hspace{4pt}        & 0.003513\hspace{4pt}        & 0.019015 & 9.00  & 0.018935 & 8.96 & 0.017643 & 8.28 \\
      10 & 0.000803$^\dagger$& \highlight{0.001442}& 0.001909$^\dagger$& 0.008893 & 10.07 & 0.007325 & 8.12 & 0.005489 & 5.84 \\
      12 & 0.000309$^\dagger$& \highlight{0.000447$^\dagger$}& 0.000810$^\dagger$& 0.007840 & 24.37 & 0.004005 & 11.96& 0.004005 &11.96 \\
      15 & 0.000152$^{\dagger\dagger}$ &  0.000152$^{\dagger\dagger}$               & 0.000152$^\dagger$& 0.007790 & 50.26 & 0.001230 & 7.09 & 0.001554 & 9.23 \\
      25 & 0.000000$^{\dagger\dagger}$ & 0.000000$^{\dagger\dagger}$ & 0.000000$^\dagger$& 0.004874 & ``$\infty$''    & 0.000213 & ``$\infty$''   &
      0.000751 & ``$\infty$''   \\
      27 & 0.000000$^{\dagger\dagger}$ & 0.000000$^{\dagger\dagger}$ & 0.000000$^\dagger$& 0.004411 & ``$\infty$''    & 0.000176 & ``$\infty$''   &
      0.000578 & ``$\infty$''   \\
      \hline
    \end{tabular}
    \caption{Performance for Penrose square root targets in the $d_1'$-metric (1958--2011 EU data)}
    \label{table_quality_d1_prime}
  \end{center}
\end{table}


\begin{table}[t]
  \begin{center}
    \begin{tabular}{|r|l|l|l|r@{}r|r@{}r|r@{}r|}
      \hline
       & \multicolumn{1}{c|}{$v^*\in \mathcal{S}$}  & \multicolumn{1}{c|}{$v^{**}\in \mathcal{C}$} & \multicolumn{1}{c|}{$v^{***}\in \mathcal{W}$}
       &\multicolumn{2}{c|}{50\%-heuristic} & \multicolumn{2}{c|}{$q^*$-heuristic}
      & \multicolumn{2}{c|}{$\bar{q}$-heuristic} \\
      \multicolumn{1}{|c|}{$n$} & \multicolumn{1}{c|}{$d_\infty$} & \multicolumn{1}{c|}{$d_\infty$} & \multicolumn{1}{c|}{$d_\infty$}
      & \multicolumn{1}{c}{$d_\infty$} & $\mathcal{S}$-error & \multicolumn{1}{c}{$d_\infty$} & $\mathcal{S}$-error & \multicolumn{1}{c}{$d_\infty$}
      & $\mathcal{S}$-error \\
      \hline
       6 & 0.014948\hspace{4pt}        & 0.014948\hspace{4pt}        & 0.014948\hspace{4pt}        & 0.082758 & 4.54  & 0.032728 & 1.19  & 0.032728 & 1.19  \\
       9 & 0.001498$^\dagger$& 0.001840\hspace{4pt}        & 0.002240\hspace{4pt}        & 0.019238 & 11.84 & 0.015909 & 9.62  & 0.023179 & 14.47 \\
      10 & 0.000575$^\dagger$& \highlight{0.001211} & 0.001960$^\dagger$& 0.011574 & 19.13 & 0.006316 & 9.98  & 0.009721 & 15.91 \\
      12 & 0.000229$^\dagger$& \highlight{0.000138$^\dagger$}& 0.000865$^\dagger$& 0.007940 & 33.67 & 0.005756 & 24.13 & 0.005756 & 24.13 \\
      15 & 0.000066$^{\dagger\dagger}$ & 0.000066$^{\dagger\dagger}$ & 0.000066$^\dagger$& 0.005923 & 88.74 & 0.001798 & 26.24 & 0.001202 & 17.21 \\
      25 & 0.000000$^{\dagger\dagger}$ & 0.000000$^{\dagger\dagger}$ & 0.000000$^\dagger$& 0.003834 & ``$\infty$''    & 0.000173 & ``$\infty$''    &
      0.000384 & ``$\infty$''    \\
      27 & 0.000000$^{\dagger\dagger}$ & 0.000000$^{\dagger\dagger}$ & 0.000000$^\dagger$& 0.003434 & ``$\infty$''    & 0.000156 & ``$\infty$''    &
      0.000277 & ``$\infty$''    \\
      \hline
    \end{tabular}
    \caption{Performance for Penrose square root targets in the $d_\infty$-metric (1958--2011 EU data)}
    \label{table_quality_d2}
  \end{center}
\end{table}

\subsection{Finite grid of objective vectors}\label{subsec:grid}


Every vector in $\mathbb{R}_{\ge 0}^n$ whose entries sum to $1$, \highlight{i.e., each element of $\Delta(n-1)$,} can in principle
be a desired power distribution in a specific context.
We approximate this infinite space by a \highlight{finite set $G^n$}. We impose $\beta_1\geq \beta_2\geq \ldots \geq \beta_n$
and let the desired power of the first $n-1$ voters be an integral multiple of $s=0.01$; the desired power of the $n$-th voter
\highlight{follows from}
the sum condition. \highlight{We refer to $G^n$ as our \emph{grid} of target vectors and to each $\beta\in G^n$ as a
\emph{grid point}.}\footnote{Step size $s$ has to be chosen with care: the
number of grid points can be intractably great already for small $n$ if $s$ is too small. But a larger $s$ induces a coarser grid of feasible
target vectors. This becomes more and more problematic as $n$ increases because of the corresponding natural decrease of an individual voter's
relative power (on average equal to ${1}/{n}$). \highlight{Choosing $s=0.25$, for instance, would result in the four  grid points $(0.5,0.25,0.25)$,
$(0.5,0.5,0)$, $(0.75,0.25,0)$, and $(1,0,0)$ for $n=3$. And $G^n$ would contain merely five grid points for \emph{any} $n\ge 4$: $(0.25,0.25,0.25,0.25,0,\dots,0)$,
$(0.5,0.25,0.25,0,\dots,0)$, $(0.5,0.5,0,\dots,0)$, $(0.75,0.25,0,\dots,0)$, and  $(1,0,\dots,0)$.}} Table~\ref{table_grid_optimal} reports key
statistics for the distribution of unavoidable deviations
from the ideal vectors in the $d_1$ and $d_\infty$-metrics: its median, average, 10\%, 5\%, and 1\%-percentile. The deviation figures are
based on the \highlight{enumerated} exact solutions in $\mathcal{W}$ for $n\le 7$ and approximations thereof for larger $n$. \highlight{For instance, the number
\highlight{0.01077} for $n=5$ in the right-most column indicates that for 1\% of the considered \highlight{46262} different target vectors $\beta$
one can obtain a distance $d_\infty(\beta, B(v^{***}))\le 0.011$, and the remaining 99\% target vectors can only be approximated less well within
$\mathcal{W}$.} A number of grid points in parentheses indicates the size of the considered random sample whenever only a subset of all grid points
could be dealt with. The deviation statistics in the corresponding rows (in light color) involve a sample error in addition to the small error of
using a conventional local \highlight{hill-climbing} algorithm instead of global optimization in $\mathcal{W}$.
\highlight{For example, the number \highlight{0.0011} for $n=15$ in the third column indicates that \emph{half} of the 10000 target vectors
$\beta$ which were sampled at random (with replacement) could, by some weighted voting game, be achieved with a $d_1$-distance of \highlight{0.0011}
or less; the remaining draws resulted in target vectors for which our search algorithm terminated with a best achievable $B(v^{***})$ further away.}

\begin{table}[htp]
  \begin{center}
    \begin{tabular}{|r|r|rrrrr|rrrrr|}
      \hline
          & \multicolumn{1}{c|}{\#grid}   & \multicolumn{5}{c|}{$d_1$-metric} & \multicolumn{5}{c|}{$d_\infty$-metric}\\
      $n$ & \multicolumn{1}{c|}{points} & \multicolumn{1}{c}{med.}   & \multicolumn{1}{c}{av.}    & \multicolumn{1}{c}{10\%} & \multicolumn{1}{c}{5\%}
      & \multicolumn{1}{c|}{1\%}  & \multicolumn{1}{c}{med.} & \multicolumn{1}{c}{av.} & \multicolumn{1}{c}{10\%} & \multicolumn{1}{c}{5\%} & \multicolumn{1}{c|}{1\%} \\
      \hline
      \!\!2\!\!&\!\!51\!\!&\!0.2400\!&\!0.2451\!&\!0.0400\!&\!0.0200\!&\!0.0000\!\!&\!\!0.12000\!&\!0.12255\!&\!0.02000\!&\!0.01000\!&\!0.00000\!\!\\
       \!\!3\!\!&\!\!884\!\!&\!0.2400\!&\!0.2278\!&\!0.1000\!&\!0.0667\!&\!0.0200\!\!&\!\!0.12000\!&\!0.11391\!&\!0.05000\!&\!0.03333\!&\!0.01000\!\!\\
       \!\!4\!\!&\!\!8037\!\!&\!0.1600\!&\!0.1622\!&\!0.0800\!&\!0.0600\!&\!0.0400\!\!&\!\!0.07000\!&\!0.07131\!&\!0.03667\!&\!0.03000\!&\!0.01500\!\!\\
       \!\!5\!\!&\!\!46262\!\!&\!0.1010\!&\!0.1135\!&\!0.0600\!&\!0.0509\!&\!0.0324\!\!&\!\!0.04000\!&\!0.04292\!&\!0.02273\!&\!0.02000\!&\!0.01077\!\!\\
       \!\!6\!\!&\!\!189509\!\!&\!0.0667\!&\!0.0790\!&\!0.0400\!&\!0.0356\!&\!0.0200\!\!&\!\!0.02222\!&\!0.02630\!&\!0.01333\!&\!0.01069\!&\!0.00815\!\!\\
       \!\!7\!\!&\!\!596763\!\!&\!0.0422\!&\!0.0543\!&\!0.0257\!&\!0.0213\!&\!0.0165\!\!&\!\!0.01255\!&\!0.01629\!&\!0.00762\!&\!0.00667\!&\!0.00495\!\!\\
       \!\!8\!\!&\!\!(10000)\!\!&\!\!\mygrey{0.0226}\!&\!\mygrey{0.0248}\!&\!\mygrey{0.0154}\!&\!\mygrey{0.0137}\!&\!\mygrey{0.0108}\!\!&\!\!\mygrey{0.00601}\!&\!\mygrey{0.00661}\!&\!\mygrey{0.00404}\!&\!\mygrey{0.00358}\!&\!\mygrey{0.00281}\!\!\\
       \!\!9\!\!&\!\!(10000)\!\!&\!\!\mygrey{0.0148}\!&\!\mygrey{0.0161}\!&\!\mygrey{0.0100}\!&\!\mygrey{0.0089}\!&\!\mygrey{0.0070}\!\!&\!\!\mygrey{0.00357}\!&\!\mygrey{0.00393}\!&\!\mygrey{0.00241}\!&\!\mygrey{0.00216}\!&\!\mygrey{0.00169}\!\!\\
      \!\!10\!\!&\!\!(10000)\!\!&\!\!\mygrey{0.0097}\!&\!\mygrey{0.0107}\!&\!\mygrey{0.0065}\!&\!\mygrey{0.0059}\!&\!\mygrey{0.0046}\!\!&\!\!\mygrey{0.00216}\!&\!\mygrey{0.00239}\!&\!\mygrey{0.00145}\!&\!\mygrey{0.00129}\!&\!\mygrey{0.00103}\!\!\\
      \!\!11\!\!&\!\!(10000)\!\!&\!\!\mygrey{0.0064}\!&\!\mygrey{0.0070}\!&\!\mygrey{0.0043}\!&\!\mygrey{0.0038}\!&\!\mygrey{0.0031}\!\!&\!\!\mygrey{0.00131}\!&\!\mygrey{0.00146}\!&\!\mygrey{0.00088}\!&\!\mygrey{0.00079}\!&\!\mygrey{0.00064}\!\!\\
      \!\!12\!\!&\!\!(10000)\!\!&\!\!\mygrey{0.0041}\!&\!\mygrey{0.0045}\!&\!\mygrey{0.0028}\!&\!\mygrey{0.0024}\!&\!\mygrey{0.0019}\!\!&\!\!\mygrey{0.00079}\!&\!\mygrey{0.00088}\!&\!\mygrey{0.00052}\!&\!\mygrey{0.00047}\!&\!\mygrey{0.00037}\!\!\\
      \!\!13\!\!&\!\!(10000)\!\!&\!\!\mygrey{0.0026}\!&\!\mygrey{0.0029}\!&\!\mygrey{0.0017}\!&\!\mygrey{0.0016}\!&\!\mygrey{0.0013}\!\!&\!\!\mygrey{0.00047}\!&\!\mygrey{0.00053}\!&\!\mygrey{0.00032}\!&\!\mygrey{0.00028}\!&\!\mygrey{0.00023}\!\!\\
      \!\!14\!\!&\!\!(10000)\!\!&\!\!\mygrey{0.0016}\!&\!\mygrey{0.0018}\!&\!\mygrey{0.0011}\!&\!\mygrey{0.0010}\!&\!\mygrey{0.0008}\!\!&\!\!\mygrey{0.00028}\!&\!\mygrey{0.00032}\!&\!\mygrey{0.00019}\!&\!\mygrey{0.00017}\!&\!\mygrey{0.00014}\!\!\\
      \!\!15\!\!&\!\!(10000)\!\!&\!\!\mygrey{0.0011}\!&\!\mygrey{0.0012}\!&\!\mygrey{0.0008}\!&\!\mygrey{0.0007}\!&\!\mygrey{0.0006}\!\!&\!\!\mygrey{0.00017}\!&\!\mygrey{0.00019}\!&\!\mygrey{0.00012}\!&\!\mygrey{0.00011}\!&\!\mygrey{0.00009}\!\!\\
      \!\!\highlight{16}\!\!&\!\!(10000)\!\!&\!\!\mygrey{0.0007}\!&\!\mygrey{0.0008}\!&\!\mygrey{0.0005}\!&\!\mygrey{0.0005}\!&\!\mygrey{0.0004}\!\!&\!\!\mygrey{0.00011}\!&\!\mygrey{0.00012}\!&\!\mygrey{0.00009}\!&\!\mygrey{0.00008}\!&\!\mygrey{0.00007}\!\!\\
      \!\!\highlight{17}\!\!&\!\!(10000)\!\!&\!\!\mygrey{0.0006}\!&\!\mygrey{0.0006}\!&\!\mygrey{0.0004}\!&\!\mygrey{0.0004}\!&\!\mygrey{0.0003}\!\!&\!\!\mygrey{0.00009}\!&\!\mygrey{0.00009}\!&\!\mygrey{0.00007}\!&\!\mygrey{0.00007}\!&\!\mygrey{0.00006}\!\!\\
      \!\!\highlight{18}\!\!&\!\!(10000)\!\!&\!\!\mygrey{0.0005}\!&\!\mygrey{0.0005}\!&\!\mygrey{0.0004}\!&\!\mygrey{0.0004}\!&\!\mygrey{0.0003}\!\!&\!\!\mygrey{0.00008}\!&\!\mygrey{0.00008}\!&\!\mygrey{0.00006}\!&\!\mygrey{0.00006}\!&\!\mygrey{0.00005}\!\!\\
      \!\!\highlight{19}\!\!&\!\!(10000)\!\!&\!\!\mygrey{0.0005}\!&\!\mygrey{0.0005}\!&\!\mygrey{0.0004}\!&\!\mygrey{0.0004}\!&\!\mygrey{0.0003}\!\!&\!\!\mygrey{0.00008}\!&\!\mygrey{0.00008}\!&\!\mygrey{0.00006}\!&\!\mygrey{0.00006}\!&\!\mygrey{0.00005}\!\!\\
      \!\!\highlight{20}\!\!&\!\!(10000)\!\!&\!\!\mygrey{0.0005}\!&\!\mygrey{0.0005}\!&\!\mygrey{0.0004}\!&\!\mygrey{0.0004}\!&\!\mygrey{0.0003}\!\!&\!\!\mygrey{0.00007}\!&\!\mygrey{0.00007}\!&\!\mygrey{0.00006}\!&\!\mygrey{0.00005}\!&\!\mygrey{0.00004}\!\!\\
      \hline
    \end{tabular}
    \caption{Distribution of unavoidable absolute deviations \highlight{$d_1(\beta,B(v^{***}))$ and $d_\infty(\beta,B(v^{***}))$}}
    \label{table_grid_optimal}
  \end{center}
\end{table}

Tables~\ref{table_grid_q_circ}--\ref{table_grid_q_bar} report analogous statistics for the distribution of absolute distances for the
three heuristics (considering each grid point for \highlight{up to $n=20$}).\footnote{\highlight{The computation times behind the unavoidable errors in Table~\ref{table_grid_optimal} ranged from less than 1\,s or 40\,m for the $n=5$ and $n=7$ rows, respectively, to 14\,h for the $n=8$ row and 31\,h for $n=20$.
Quite some time is spent on approximating the exact solution of the inverse problem. Times for the heuristic in, e.g., Table~\ref{table_grid_q_star}
were only 1\,s, 15\,s, 50\,m, and 3\,h, respectively. A sample of 10000 grid points represents a reasonable compromise between precision and computational effort. Raising the sample size to 100000 would, e.g, have produced the median, average and quantile entries $(0.0041, 0.0045, 0.0027, 0.0025, 0.0020)$ for the $d_1$-metric and $n=12$; lowering it to 1000 would have resulted in $(0.0042, 0.0047, 0.0028, 0.0026, 0.0018)$.}} A comparison of the respective
deviation statistics with those in Table~\ref{table_grid_optimal} broadly confirm the observations that were made for the very specific \highlight{target vectors}
derived from Penrose's square root rule in Section~\ref{subsec:EU}: the average and each reported percentile of the avoidable deviations
decrease in $n$.
They can be regarded as small in absolute terms, but they are sizeable in relative terms. Again the 50\%-heuristic
is clearly outperformed (in the sense of first order stochastic dominance) by the $q^*$ and $\bar q$-heuristics for $n\ge 3$.

\begin{table}[t]
  \begin{center}
    \begin{tabular}{|r|r|rrrrr|rrrrr|}
      \hline
          & \multicolumn{1}{c|}{\#grid}   & \multicolumn{5}{c|}{$d_1$-metric} & \multicolumn{5}{c|}{$d_\infty$-metric}\\
      $n$ & \multicolumn{1}{c|}{points} & \multicolumn{1}{c}{med.}   & \multicolumn{1}{c}{av.}    & \multicolumn{1}{c}{10\%} &
      \multicolumn{1}{c}{5\%}  & \multicolumn{1}{c|}{1\%}  & \multicolumn{1}{c}{med.}
      & \multicolumn{1}{c}{av.} & \multicolumn{1}{c}{10\%} & \multicolumn{1}{c}{5\%} & \multicolumn{1}{c|}{1\%} \\
      \hline
       2 &        51 & 0.480 & 0.480 & 0.080 & 0.020 & 0.000 & 0.240 & 0.240 & 0.040 & 0.010 & 0.000 \\
       3 &       884 & 0.560 & 0.555 & 0.200 & 0.133 & 0.047 & 0.280 & 0.278 & 0.100 & 0.067 & 0.023 \\
       4 &      8037 & 0.440 & 0.509 & 0.200 & 0.153 & 0.080 & 0.210 & 0.249 & 0.083 & 0.063 & 0.033 \\
       5 &     46262 & 0.347 & 0.448 & 0.160 & 0.127 & 0.075 & 0.153 & 0.209 & 0.061 & 0.049 & 0.029 \\
       6 &    189509 & 0.297 & 0.389 & 0.129 & 0.103 & 0.066 & 0.120 & 0.177 & 0.045 & 0.035 & 0.023 \\
       7 &    596763 & 0.247 & 0.338 & 0.101 & 0.080 & 0.052 & 0.097 & 0.151 & 0.033 & 0.025 & 0.016 \\
       8 &   1527675 & 0.206 & 0.297 & 0.080 & 0.063 & 0.041 & 0.080 & 0.132 & 0.025 & 0.019 & 0.012 \\
       9 &   3314203 & 0.176 & 0.265 & 0.064 & 0.051 & 0.034 & 0.068 & 0.118 & 0.020 & 0.015 & 0.009 \\
      10 &   6292069 & 0.153 & 0.240 & 0.053 & 0.043 & 0.029 & 0.059 & 0.107 & 0.016 & 0.012 & 0.007 \\
      11 &  10718685 & 0.136 & 0.220 & 0.046 & 0.037 & 0.025 & 0.052 & 0.099 & 0.014 & 0.010 & 0.006 \\
      12 &  \highlight{16713148} & 0.123 & 0.205 & 0.041 & 0.033 & 0.023 & 0.047 & 0.092 & 0.012 & 0.009 & 0.005 \\
      13 &  \highlight{24234058} & 0.112 & 0.193 & 0.038 & 0.030 & 0.021 & 0.044 & 0.087 & 0.011 & 0.008 & 0.005 \\
      14 &  \highlight{33097743} & 0.104 & 0.183 & 0.035 & 0.028 & 0.020 & 0.041 & 0.083 & 0.010 & 0.008 & 0.004 \\
      15 &  \highlight{43018955} & 0.097 & 0.175 & 0.033 & 0.027 & 0.019 & 0.038 & 0.079 & 0.010 & 0.007 & 0.004 \\
\highlight{16} & \highlight{53662038}& \highlight{0.092} & \highlight{0.169} & \highlight{0.032} & \highlight{0.026} & \highlight{0.018} & \highlight{0.037} & \highlight{0.076} & \highlight{0.009} & \highlight{0.007} & \highlight{0.004} \\
\highlight{17} & \highlight{64684584}& \highlight{0.087} & \highlight{0.164} & \highlight{0.031} & \highlight{0.025} & \highlight{0.017} & \highlight{0.035} & \highlight{0.074} & \highlight{0.009} & \highlight{0.006} & \highlight{0.004} \\
\highlight{18} & \highlight{75772412}& \highlight{0.084} & \highlight{0.159} & \highlight{0.030} & \highlight{0.024} & \highlight{0.017} & \highlight{0.034} & \highlight{0.072} & \highlight{0.009} & \highlight{0.006} & \highlight{0.003} \\
\highlight{19} & \highlight{86658411}& \highlight{0.081} & \highlight{0.156} & \highlight{0.029} & \highlight{0.024} & \highlight{0.016} & \highlight{0.033} & \highlight{0.071} & \highlight{0.008} & \highlight{0.006} & \highlight{0.003} \\
\highlight{20} & \highlight{97132873}& \highlight{0.078} & \highlight{0.153} & \highlight{0.028} & \highlight{0.023} & \highlight{0.016} & \highlight{0.032} & \highlight{0.070} & \highlight{0.008} & \highlight{0.006} & \highlight{0.003} \\
      \hline
    \end{tabular}
    \caption{Distribution of absolute deviations for the 50\%-heuristic}
    \label{table_grid_q_circ}
  \end{center}
\end{table}

\begin{table}[t]
  \begin{center}
    \begin{tabular}{|r|r|rrrrr|rrrrr|}
      \hline
          & \multicolumn{1}{c|}{\#grid}   & \multicolumn{5}{c|}{$d_1$-metric} & \multicolumn{5}{c|}{$d_\infty$-metric}\\
      $n$ & \multicolumn{1}{c|}{points} & \multicolumn{1}{c}{med.}   & \multicolumn{1}{c}{av.}    & \multicolumn{1}{c}{10\%} &
      \multicolumn{1}{c}{5\%}  & \multicolumn{1}{c|}{1\%}  &
      \multicolumn{1}{c}{med.}   & \multicolumn{1}{c}{av.} & \multicolumn{1}{c}{10\%} & \multicolumn{1}{c}{5\%} & \multicolumn{1}{c|}{1\%} \\     \hline
       2 &        51 & 0.480 & 0.480 & 0.080 & 0.020 & 0.000 & 0.240 & 0.240 & 0.040 & 0.010 & 0.000 \\
       3 &       884 & 0.400 & 0.434 & 0.160 & 0.107 & 0.040 & 0.200 & 0.217 & 0.080 & 0.053 & 0.020 \\
       4 &      8037 & 0.340 & 0.370 & 0.160 & 0.120 & 0.060 & 0.147 & 0.172 & 0.065 & 0.050 & 0.025 \\
       5 &     46262 & 0.280 & 0.312 & 0.133 & 0.107 & 0.062 & 0.113 & 0.138 & 0.052 & 0.040 & 0.023 \\
       6 &    189509 & 0.227 & 0.263 & 0.109 & 0.088 & 0.058 & 0.088 & 0.112 & 0.038 & 0.030 & 0.020 \\
       7 &    596763 & 0.189 & 0.224 & 0.087 & 0.070 & 0.047 & 0.071 & 0.093 & 0.027 & 0.021 & 0.014 \\
       8 &   1527675 & 0.158 & 0.192 & 0.066 & 0.053 & 0.035 & 0.056 & 0.079 & 0.020 & 0.015 & 0.010 \\
       9 &   3314203 & 0.133 & 0.168 & 0.051 & 0.040 & 0.026 & 0.047 & 0.068 & 0.014 & 0.011 & 0.007 \\
      10 &   6292069 & 0.114 & 0.148 & 0.039 & 0.030 & 0.019 & 0.039 & 0.060 & 0.011 & 0.008 & 0.005 \\
      11 &  10718685 & 0.098 & 0.132 & 0.030 & 0.023 & 0.014 & 0.033 & 0.054 & 0.008 & 0.006 & 0.003 \\
      12 &  \highlight{16713148} & 0.086 & 0.120 & 0.024 & 0.017 & 0.010 & 0.029 & 0.049 & 0.006 & 0.004 & 0.002 \\
      13 &  \highlight{24234058} & 0.075 & 0.110 & 0.019 & 0.013 & 0.007 & 0.026 & 0.045 & 0.005 & 0.003 & 0.002 \\
      14 &  \highlight{33097743} & 0.068 & 0.102 & 0.016 & 0.011 & 0.005 & 0.023 & 0.042 & 0.004 & 0.003 & 0.001 \\
      15 &  \highlight{43018955} & 0.061 & 0.096 & 0.013 & 0.008 & 0.004 & 0.021 & 0.040 & 0.003 & 0.002 & 0.001 \\
\highlight{16} & \highlight{53662038}& \highlight{0.056} & \highlight{0.091} & \highlight{0.011} & \highlight{0.007} & \highlight{0.003} & \highlight{0.019} & \highlight{0.038} & \highlight{0.003} & \highlight{0.002} & \highlight{0.001} \\
\highlight{17} & \highlight{64684584}& \highlight{0.052} & \highlight{0.087} & \highlight{0.009} & \highlight{0.006} & \highlight{0.003} & \highlight{0.018} & \highlight{0.036} & \highlight{0.002} & \highlight{0.001} & \highlight{0.001} \\
\highlight{18} & \highlight{75772412}& \highlight{0.049} & \highlight{0.083} & \highlight{0.008} & \highlight{0.005} & \highlight{0.002} & \highlight{0.017} & \highlight{0.035} & \highlight{0.002} & \highlight{0.001} & \highlight{0.000} \\
\highlight{19} & \highlight{86658411}& \highlight{0.046} & \highlight{0.081} & \highlight{0.007} & \highlight{0.004} & \highlight{0.002} & \highlight{0.016} & \highlight{0.034} & \highlight{0.002} & \highlight{0.001} & \highlight{0.000} \\
\highlight{20} & \highlight{97132873}& \highlight{0.044} & \highlight{0.078} & \highlight{0.006} & \highlight{0.004} & \highlight{0.002} & \highlight{0.015} & \highlight{0.033} & \highlight{0.002} & \highlight{0.001} & \highlight{0.000} \\
      \hline
    \end{tabular}
    \caption{Distribution of absolute deviations for the $q^\star$-heuristic}
    \label{table_grid_q_star}
  \end{center}
\end{table}

\begin{table}[t]
  \begin{center}
    \begin{tabular}{|r|r|rrrrr|rrrrr|}
      \hline
          & \multicolumn{1}{c|}{\#grid}   & \multicolumn{5}{c|}{$d_1$-metric} & \multicolumn{5}{c|}{$d_\infty$-metric}\\
      $n$ & \multicolumn{1}{c|}{points} & \multicolumn{1}{c}{med.}   & \multicolumn{1}{c}{av.}    & \multicolumn{1}{c}{10\%} &
      \multicolumn{1}{c}{5\%}  & \multicolumn{1}{c|}{1\%}  &
      \multicolumn{1}{c}{med.}   & \multicolumn{1}{c}{av.} & \multicolumn{1}{c}{10\%} & \multicolumn{1}{c}{5\%} & \multicolumn{1}{c|}{1\%} \\
      \hline
       2 &        51 & 0.280 & 0.327 & 0.040 & 0.020 & 0.000 & 0.140 & 0.164 & 0.020 & 0.010 & 0.000 \\
       3 &       884 & 0.320 & 0.332 & 0.140 & 0.100 & 0.040 & 0.160 & 0.166 & 0.070 & 0.050 & 0.020 \\
       4 &      8037 & 0.300 & 0.304 & 0.147 & 0.110 & 0.050 & 0.130 & 0.138 & 0.063 & 0.045 & 0.020 \\
       5 &     46262 & 0.250 & 0.263 & 0.132 & 0.100 & 0.060 & 0.101 & 0.111 & 0.050 & 0.040 & 0.023 \\
       6 &    189509 & 0.204 & 0.224 & 0.104 & 0.085 & 0.056 & 0.077 & 0.088 & 0.036 & 0.029 & 0.019 \\
       7 &    596763 & 0.159 & 0.187 & 0.079 & 0.064 & 0.043 & 0.057 & 0.069 & 0.024 & 0.020 & 0.013 \\
       8 &   1527675 & 0.120 & 0.155 & 0.058 & 0.047 & 0.033 & 0.041 & 0.055 & 0.017 & 0.013 & 0.009 \\
       9 &   3314203 & 0.098 & 0.133 & 0.044 & 0.036 & 0.024 & 0.033 & 0.046 & 0.012 & 0.010 & 0.006 \\
      10 &   6292069 & 0.079 & 0.115 & 0.033 & 0.026 & 0.018 & 0.026 & 0.039 & 0.009 & 0.007 & 0.004 \\
      11 &  10718685 & 0.071 & 0.103 & 0.027 & 0.021 & 0.014 & 0.023 & 0.035 & 0.007 & 0.005 & 0.003 \\
      12 &  \highlight{16713148} & 0.058 & 0.092 & 0.021 & 0.016 & 0.010 & 0.018 & 0.030 & 0.005 & 0.004 & 0.002 \\
      13 &  \highlight{24234058} & 0.050 & 0.084 & 0.016 & 0.012 & 0.007 & 0.015 & 0.027 & 0.004 & 0.003 & 0.002 \\
      14 &  \highlight{33097743} & 0.045 & 0.078 & 0.014 & 0.010 & 0.006 & 0.014 & 0.025 & 0.003 & 0.002 & 0.001 \\
      15 &  \highlight{43018955} & 0.042 & 0.074 & 0.012 & 0.009 & 0.005 & 0.012 & 0.024 & 0.003 & 0.002 & 0.001 \\
\highlight{16} & \highlight{53662038}& \highlight{0.039} & \highlight{0.070} & \highlight{0.011} & \highlight{0.007} & \highlight{0.004} & \highlight{0.011} & \highlight{0.022} & \highlight{0.002} & \highlight{0.002} & \highlight{0.001} \\
\highlight{17} & \highlight{64684584}& \highlight{0.040} & \highlight{0.070} & \highlight{0.011} & \highlight{0.008} & \highlight{0.004} & \highlight{0.012} & \highlight{0.023} & \highlight{0.002} & \highlight{0.002} & \highlight{0.001} \\
\highlight{18} & \highlight{75772412}& \highlight{0.037} & \highlight{0.067} & \highlight{0.009} & \highlight{0.007} & \highlight{0.003} & \highlight{0.011} & \highlight{0.022} & \highlight{0.002} & \highlight{0.001} & \highlight{0.001} \\
\highlight{19} & \highlight{86658411}& \highlight{0.041} & \highlight{0.069} & \highlight{0.011} & \highlight{0.008} & \highlight{0.004} & \highlight{0.013} & \highlight{0.024} & \highlight{0.002} & \highlight{0.002} & \highlight{0.001} \\
\highlight{20} & \highlight{97132873}& \highlight{0.038} & \highlight{0.067} & \highlight{0.010} & \highlight{0.007} & \highlight{0.003} & \highlight{0.012} & \highlight{0.023} & \highlight{0.002} & \highlight{0.001} & \highlight{0.001} \\
      \hline
    \end{tabular}
    \caption{Distribution of absolute deviations for the $\bar{q}$-heuristic}
    \label{table_grid_q_bar}
  \end{center}
\end{table}

\subsection{Analytical example}\label{subsec:analytical}
\highlight{A statement which would be analogous to observation~(v) in Section~\ref{subsec:EU} cannot be deduced on the basis of statistical information as provided by Tables~\ref{table_grid_q_circ}--\ref{table_grid_q_bar}. We, therefore, close our computational investigation by studying a particularly simple analytical example. It shows  transparently that~-- as indicated by observation~(v)~-- \emph{relative} deviations between the considered heuristics and globally optimal solutions need not disappear for
$n\to \infty$.}

\highlight{Consider the desired power distribution}
$$
  \beta^n=\frac{1}{2n-1}(\underbrace{2,\dots,2}_{n-1\textnormal{ twos}},1)
$$
for $n\geq 2$ \highlight{and choose $w^n=\beta^n$ as all three heuristics do}.\footnote{The construction is inspired by a sequence of weighted voting games to which the Penrose limit theorem does
\emph{not} apply \highlight{even though every voter's relative weight vanishes as $n\to \infty$. Namely, the sequence $\{(\frac{1}{2};w^n)\}_{n\in \mathbb{N}}$ belongs to the class of somewhat pathological examples alluded to on p.~\pageref{rem:pathological} (cf.\ \citeNP{Lindner/Owen:2007})}.} For any quota $q\in I_1^j=\frac{1}{2n-1}\cdot (2j-1,2j]$, where $1\le j\le n-1$ and $j\in\mathbb{N}$, the PBI of the
smallest constituency is exactly zero and, by symmetry, the \highlight{(normalized)} PBI of \highlight{each of} the other constituencies equals $\frac{1}{n-1}$.
For the remaining possibilities $q\in I_2^j=\frac{1}{2n-1}\cdot (2j,2j+1]$ where $0\le j\le n-1$, all constituencies
have a PBI of $\frac{1}{n}$. Denoting the corresponding weighted games by $v_{1,j}^n$ and $v_{2,j}^n$ one obtains
\begin{eqnarray*}
  && d_1\!\left(v_{1,j}^n,\beta^n\right)=\frac{2}{2n-1},\\
  && d_1\!\left(v_{2,j}^n,\beta^n\right)=\frac{2}{2n-1}\cdot\frac{n-1}{n},\\
  && d_\infty\!\left(v_{1,j}^n,\beta^n\right)=\frac{1}{2n-1},\textnormal{ and}\\
  && d_\infty\!\left(v_{2,j}^n,\beta^n\right)=\frac{1}{2n-1}\cdot\frac{n-1}{n}.
\end{eqnarray*}
So, independently of the quota, the $\Vert\cdot\Vert_1$-error is $\frac{2}{2n-1}+O\!\left(n^{-2}\right)$ and the
$\Vert\cdot\Vert_\infty$-error is $\frac{1}{2n-1}+O\!\left(n^{-2}\right)$.

The $\bar{q}$ and $q^*$-heuristics prescribe quotas of
\begin{eqnarray*}
  && \bar{q}=\frac{1}{2}+\frac{1}{\sqrt{\pi n}},\textnormal{ and}\\
  && q^*=\frac{1}{2}+\frac{\sqrt{4n-3}}{4n-2},
\end{eqnarray*}
respectively. They and $q^\circ=50\%$ fall into $I_1^j$ and $I_2^j$ for infinitely many $n$. Thus, all three rules render
the smallest constituency a null player infinitely many times as $n\to \infty$, just as it happened to Luxembourg in the
EEC Council between 1958 and 1973, \highlight{and yield the indicated distances}.

\begin{table}[t]
  \begin{center}
    \begin{tabular}{|r|l|l|r|rr|}
      \hline
       &\multicolumn{1}{c|}{$v^*\in \mathcal{S}$}  & \multicolumn{1}{c|}{$v^{**}\in \mathcal{C}$} & \multicolumn{1}{c|}{$v^{***}\in \mathcal{W}$}
       &\multicolumn{2}{c|}{$q$-heuristic}\\
      \multicolumn{1}{|c|}{$n$} & \multicolumn{1}{c|}{$d_1$}  & \multicolumn{1}{c|}{$d_1$} & \multicolumn{1}{c|}{$d_1$} & \multicolumn{1}{c}{$d_1$}
      & \multicolumn{1}{c|}{$\mathcal{C}$-error} \\
      \hline
       2 & 0.333333       & 0.333333\hspace{4pt} & 0.333333\hspace{4pt} & 0.333333        & 0.000000     \\
       3 & 0.266667       & 0.266667\hspace{4pt} & 0.266667\hspace{4pt} & 0.266667        & 0.000000     \\
       4 & 0.214286       & 0.214286\hspace{4pt} & 0.214286\hspace{4pt} & 0.214286        & 0.000000     \\
       5 & 0.038647       & 0.158730\hspace{4pt} & 0.158730\hspace{4pt} & 0.177778        & 0.120000     \\
       6 & 0.000000       & 0.113636\hspace{4pt} & 0.113636\hspace{4pt} & 0.151515        & 0.333333     \\
       7 & 0.000000       & 0.085470\hspace{4pt} & 0.085470\hspace{4pt} & 0.131868        & 0.542857     \\
       8 & 0.000000       & 0.066667\hspace{4pt} & 0.066667\hspace{4pt} & 0.116667        & 0.750000     \\
       9 & 0.000000       & 0.064171\hspace{4pt} & 0.064171\hspace{4pt} & 0.104575        & 0.629630     \\
      10 & 0.000000       & 0.061042\hspace{4pt} & 0.061042\hspace{4pt} & 0.094737        & 0.552000     \\
      11 & 0.000000       & 0.052158\hspace{4pt} & 0.052158\hspace{4pt} & 0.086580        & 0.659944     \\
      12 & 0.000000       & 0.047254\hspace{4pt} & 0.047254\hspace{4pt} & 0.079710        & 0.686856     \\
      13 & 0.000000       & \highlight{0.042353}\hspace{4pt} & \highlight{0.042353\hspace{4pt}}  & 0.073846        & \highlight{0.743590}     \\
\highlight{14} & \highlight{0.000000} & \highlight{0.037037\hspace{4pt}}  & \highlight{0.037037\hspace{4pt}} & \highlight{0.068783} & \highlight{0.857143} \\
\highlight{15} & \highlight{0.000000} & \highlight{0.034483$^{\dagger\dagger}$} & \highlight{0.034483$^\dagger$} & \highlight{0.064368}        & \highlight{0.866667} \\
\highlight{16} & \highlight{0.000000} & \highlight{0.033845$^{\dagger\dagger}$} & \highlight{0.033845$^\dagger$} & \highlight{0.060484}        & \highlight{0.780576} \\
\highlight{17} & \highlight{0.000000} & \highlight{0.032221$^{\dagger\dagger}$} & \highlight{0.032221$^\dagger$} & \highlight{0.057041} & \highlight{0.770270} \\
\highlight{18} & \highlight{0.000000} & \highlight{0.030866$^{\dagger\dagger}$} & \highlight{0.030866$^\dagger$} & \highlight{0.053968} & \highlight{0.748490} \\
\highlight{19} &     & \highlight{0.028108$^{\dagger\dagger}$} & \highlight{0.028108$^\dagger$} & \highlight{0.051209} & \highlight{0.821862} \\
\highlight{20} &     & \highlight{0.025641$^{\dagger\dagger}$} & \highlight{0.025641$^\dagger$} & \highlight{0.048718} & \highlight{0.900000} \\
      \hline
    \end{tabular}
    \caption{Deviations from $\beta^n$  in the $d_1$-metric (analytical example)}
    \label{table_quality_d1_analytical}
  \end{center}
\end{table}

In contrast, there always exists a simple game $v^*\in \mathcal{S}$ whose PBI attains $\beta^n$ \emph{exactly} for
$6\le n\le \highlight{18}$.\footnote{\highlight{An example for $n=6$ is given by the following set of minimal winning
coalitions: $\{2,4,5,6\}$, $\{2,3,4,5\}$, $\{1,3,5,6\}$, $\{1,3,4,5\}$, $\{1,2,4,6\}$, and $\{1,2,3,5\}$, which attains
the PBI vector $\frac{1}{44}\left(8,8,8,8,8,4\right)=\frac{1}{11}(2,2,2,2,2,1)$.}} And we conjecture that this remains
true for $n\ge \highlight{19}$. \highlight{So the corresponding distance between $\beta$ and $B(v^*)$ equals 0 independently
of the considered metric.} Approximation results for complete simple games and distances for the heuristic choice of $w=\beta^n$ with an
``optimal'' quota $q$ that leads to $v_{2,j}^n$ (abbreviated as \emph{$q$-heuristic}) are summarized in
Tables~\ref{table_quality_d1_analytical} and \ref{table_quality_d2_analytical}.\footnote{\highlight{We remark that we have imposed suitable equivalence classes of voters in the optimizations for $n\ge 15$ in order
to reduce the computational burden. Voters $i$ and $j$ are in the same equivalence class 
if
$v(U)=v(U\cup\{j\}\backslash\{i\})$ for all coalitions $U$ with $i\in U$ and $j\notin U$. 
This is more restrictive than requiring identical PBI values of $i$ and $j$ only.}} Since the unavoidable error in the class of
simple games $\mathcal{S}$ (and hence of finite intersections of weighted games $v\in \mathcal{W}$) is zero for
$6\le n\le \highlight{18}$ and presumably beyond, we consider the $\mathcal{C}$-error in order to evaluate the relative performance of
the $q$-heuristic.

\highlight{Tables~\ref{table_quality_d1_analytical} and \ref{table_quality_d2_analytical} suggest that (i) the $\mathcal{C}$-error
converges to a positive constant in case of the $d_1$-metric and (ii) this error even seems to grow without bound when the $d_\infty$-metric is
used. The key finding that \emph{relative errors fail to disappear as $n$ grows large} can be made more rigorous. To this end, consider the
sequence of weighted voting games $\{(q^n; w^n)\}_{n\in \mathbb{N}}$ with
$$
(q^n;w^n)=(2n-a-4;\underbrace{3,\ldots, 3}_{a \textnormal{ threes}}, \underbrace{2,\ldots, 2}_{n-a-1 \textnormal{ twos}},1)
$$
for a suitable parameter $a$. If $a$ is chosen equal to about $\frac{6n}{7}$ (the exact number is provided in Lemma~\ref{lem:6n7specifics}) then the $d_1$-distance between $\beta^n$ and $B(q^n; w^n)$ asymptotically tends to
$\frac{1}{2n}$.\footnote{\highlight{Details on the simple but tedious computations are provided in Appendix~B.}}
The distance achieved by these specific weighted voting games provides an upper bound for the distance achieved by the respective optimal weighted game $v^{***}\in\mathcal{W}$. The $\mathcal{W}$-error of the heuristic under the $d_1$-metric is, therefore, asymptotically bounded from below by}
$$
\lim_{n\to \infty}\frac{\left|\frac{2}{2n-1}\cdot\frac{n-1}{n}-\frac{1}{2n}\right|}{\frac{1}{2n}}=1.
$$
\highlight{So the $\mathcal{W}$-error and a~fortiori also the $\mathcal{C}$-error stay at around 100\% even as $n\to \infty$ when discrepancies between $\beta^n$ and $B(v^{**})$ (or $B(v^{***})$) are evaluated by the $d_1$-metric.}

\highlight{Analogously, one can choose $a$ equal to about $\frac{2n}{3}$ and then check that the $d_\infty$-distance between the game's PBI and $\beta^n$ asymptotically tends to $\frac{1}{n^2}$
as $n\to \infty$. This again translates into an asymptotic lower bound for the $\mathcal{W}$-error, and a~fortiori the $\mathcal{C}$-error; these relative errors go to infinity at an approximately linear speed.}

\begin{table}[t]
  \begin{center}
    \begin{tabular}{|r|r|l|l|rr|}
      \hline
      &\multicolumn{1}{c|}{$v^*\in \mathcal{S}$}  & \multicolumn{1}{c|}{$v^{**}\in \mathcal{C}$} & \multicolumn{1}{c|}{$v^{***}\in \mathcal{W}$}
      &\multicolumn{2}{c|}{$q$-heuristic}\\
      \multicolumn{1}{|c|}{$n$} & \multicolumn{1}{c|}{$d_\infty$}  & \multicolumn{1}{c|}{$d_\infty$} & \multicolumn{1}{c|}{$d_\infty$} &
      \multicolumn{1}{c}{$d_\infty$} & \multicolumn{1}{c|}{$\mathcal{C}$-error} \\
      \hline
       2 & 0.166667       & 0.166667\hspace{4pt}           & 0.166667\hspace{4pt}         & 0.166667        & 0.000000     \\
       3 & 0.133333       & 0.133333\hspace{4pt}           & 0.133333\hspace{4pt}         & 0.133333        & 0.000000     \\
       4 & 0.107143       & 0.107143\hspace{4pt}           & 0.107143\hspace{4pt}         & 0.107143        & 0.000000     \\
       5 & 0.019324       & 0.050505\hspace{4pt}           & 0.050505\hspace{4pt}         & 0.088889        & 0.760000     \\
       6 & 0.000000       & 0.034759\hspace{4pt}           & 0.034759\hspace{4pt}         & 0.075758        & 1.179487     \\
       7 & 0.000000       & 0.022624\hspace{4pt}           & 0.022624\hspace{4pt}         & 0.065934        & 1.914286     \\
       8 & 0.000000       & 0.015686\hspace{4pt}           & 0.015686\hspace{4pt}         & 0.058333        & 2.718750     \\
       9 & 0.000000       & 0.014199\hspace{4pt}           & 0.014199\hspace{4pt}         & 0.052288        & 2.682540     \\
      10 & 0.000000       & 0.008772\hspace{4pt}           & 0.008772\hspace{4pt}         & 0.047368        & 4.400000     \\
      11 & 0.000000       & 0.008282\hspace{4pt}           & 0.008282\hspace{4pt}         & 0.043290        & 4.227273     \\
      12 & 0.000000       & 0.007688\hspace{4pt}           & 0.007688\hspace{4pt}         & 0.039855        & 4.183908     \\
      13 & 0.000000       & \highlight{0.005373\hspace{4pt}} & \highlight{0.005373\hspace{4pt}}   & 0.036923 & \highlight{5.871795} \\
      \highlight{14} & \highlight{0.000000} & \highlight{0.005109\hspace{4pt}} & \highlight{0.005109\hspace{4pt}} & \highlight{0.034392}  & \highlight{5.732143} \\
      \highlight{15} & \highlight{0.000000} & \highlight{0.004628$^\dagger$} & \highlight{0.004815$^\dagger$} & \highlight{0.032184}        & \highlight{5.954839} \\
      \highlight{16} & \highlight{0.000000} & \highlight{0.003619$^{\dagger\dagger}$} & \highlight{0.003619$^\dagger$} & \highlight{0.030242}        & \highlight{7.357143} \\
      \highlight{17} & \highlight{0.000000} & \highlight{0.003463$^{\dagger\dagger}$} & \highlight{0.003463$^\dagger$} & \highlight{0.028520} & \highlight{7.235294} \\
      \highlight{18} & \highlight{0.000000} & \highlight{0.003297$^{\dagger\dagger}$} & \highlight{0.003297$^\dagger$} & \highlight{0.026984} & \highlight{7.185185} \\
      \highlight{19} &                & \highlight{0.002600$^{\dagger\dagger}$} & \highlight{0.002600$^\dagger$} & \highlight{0.025605} & \highlight{8.848225} \\
      \highlight{20} &                & \highlight{0.002502$^{\dagger\dagger}$} & \highlight{0.002502$^\dagger$} & \highlight{0.024359} & \highlight{8.737500} \\
      \hline
    \end{tabular}
    \caption{Deviations from $\beta^n$  in the $d_\infty$-metric (analytical example)}
    \label{table_quality_d2_analytical}
  \end{center}
\end{table}
\section{Conclusion}
\label{sec_conclusion}
The computations which we have reported in Section~\ref{subsec:EU} confirm that if one wants to implement the Penrose square root
rule for population data from today's European Union, the  $q^*$-heuristic of S\l{}omczy\'nski and \.Zyczkowski and, to a lesser
extent, the even simpler $\bar q$-heuristic perform very well in absolute terms. That is, the distance between a (normalized) square
root target distribution $\beta$ and the PBI $B(q^*,\beta)$ is close to zero. However, the considered heuristics can still be very
far from the globally optimal solution to the inverse problem in relative terms. This finding applies even when only weighted voting
games are allowed as feasible solutions. And it is not restricted to small voting bodies, but holds for the current number of EU
members $n=27$.

The extensive computations reported in Section~\ref{subsec:grid} confirm this observation. They provide the first systematic evaluation
of the unavoidable deviations between arbitrary target PBI power vectors and those that are actually implementable for voting bodies
with \highlight{up to $n=20$} members. Numbers such as the ones reported in Table~\ref{table_grid_optimal} can potentially be useful
in order to improve termination criteria for local search algorithms (e.g., Leech \citeyearNP{Leech:2002,Leech:2003}), which have been
used in applied studies. If, say, a locally optimal candidate solution for an inverse problem with $n=11$ voters has a $d_1$-deviation
from the desired vector $\beta$ greater than 0.0064, then Table~\ref{table_grid_optimal} indicates that the odds of further improvements
in the class of weighted voting games are 50:50 and search presumably should continue in a different part of the game space. If, however,
the deviation is smaller than 0.0031, then the odds are rather 1:99; termination might then make sense.

That desired PBI distributions which concentrate a major share of relative power amongst a few voters pose problems for the considered
heuristics is not surprising. After all, the derivation of $q^*$ by \citeN{Slomczynski/Zyczkowski:2007} involves 
a technical condition (see fn.~\ref{fn:SZcondition}) from which one can conclude that the target PBI of a single voter should approach
zero at least as fast as $1/\sqrt{n}$. It is much less obvious, however, that, first, it is not sufficient to have a target vector $\beta$
without ``outliers'' in order to obtain a heuristic solution that is good relative to the exact one and, second, the relative errors may
get larger rather than smaller as $n$ increases. This emerged from the extensive numerical computations reported in
Sections~\ref{subsec:EU}--\ref{subsec:grid} and \highlight{has been formally demonstrated} for a specific analytical example in Section~\ref{subsec:analytical}.
One might, therefore, summarize our findings as justifying and potentially even calling for case-specific optimization rather than the
application of a generally rather good heuristic~-- not only for small but even for large voting bodies.

\section*{Acknowledgements}
We thank Nicola Maaser, Friedrich Pukelsheim and two anonymous referees for constructive comments on an earlier version. Napel acknowledges the generous hospitality of the Max Planck Institute of Economics, Jena, where parts of this paper were written.

\appendix
\section*{Appendix A: ILP formulation for the inverse Penrose-Banzhaf index problem}
\label{sec_appendix_algorithm}

\highlight{Even though stating the inverse power index problem as an optimization problem is trivial (see (\ref{eq:min_problem_stated})), coming up with an implementable formulation for finding an exact solution is not. Voting systems are discrete objects, and so some kind of discrete
optimization is needed. Exhaustive enumeration (see \citeNP{Keijzer/Klos/Zhang:2010}) is limited at best to $n\le 9$ (see Table~\ref{table_number_of_voting_systems}).
A much more tractable alternative is to describe the set of feasible binary voting systems by integer
variables and to use some of the available optimization software packages.
These allow significantly larger numbers of variables when dealing with \emph{linear} rather than \emph{non-linear} (mixed) integer optimization problems. So it is unfortunate that problem (\ref{eq:min_problem_stated}) cannot directly be translated into a linear problem for the Penrose-Banzhaf index.\footnote{Interestingly, one can easily linearize the analogous inverse problem for the Shapley-Shubik power index (SSI). So even though the PBI is easier to compute than the SSI, the corresponding inverse problem is more difficult.}
The ``work-around'', which has first been suggested by \citeN{Kurz:2012} and is adopted here, is to use ILP techniques in order to merely find out whether some binary voting system $v$ exists whose PBI vector $B(v)$ is at most a specified distance $\alpha$ apart from the target $\beta$. This \emph{feasibility problem} can be solved much more easily than the underlying minimization problem. Still, one can iteratively determine the exact solution of (\ref{eq:min_problem_stated}) by varying $\alpha$.
}

\highlight{We will mostly confine our description to the case of measuring distance by the $d_1$-metric. Adaptations to the $d_{1}'$ or
$d_\infty$-metric are straightforward. They involve heterogeneous coefficients in inequality (\ref{ie_alpha}) below for $d_1'$, and neither $i$-subscripts in (\ref{ie_delta_i})--(\ref{ie_abs_2}) nor a summation in (\ref{ie_alpha}) for $d_\infty$.}

\highlight{The PBI vector $(1,0,\dots,0)$ of a dictator has at most a $d_1$-distance of $2$ from any normalized power distribution (summing to $1$). This is in fact the worst case, and the minimal achievable deviation $\alpha^*$ must lie inside the interval $[l_1,u_1]$, where $l_1=0$ and $u_1=2$. In each iteration $t=1, \ldots, T$ of the algorithm we will check whether
$\alpha=({u_{t}-l_{t}})/{2}$ is a feasible distance between target $\beta$ and the PBI values generated by the considered class of voting systems. If so, we set $u_{t+1}=\alpha$ and leave $l_{t+1}=l_t$ unchanged; otherwise we update $l_{t+1}=\alpha$ and leave $u_{t+1}=u_t$ unchanged. In each iteration the length
of the interval $[l_t,u_t]$ shrinks by a factor of $2$. Since the total number of swings in an $n$-player voting game lies
between $n$ and $m{n\choose m}<n2^n$ where $m=\left\lfloor\frac{n}{2}\right\rfloor+1$ (see, e.g., \citeNP[sec.~3.3]{Felsenthal/Machover:1998}),
two distinct PBI  vectors differ, both in the $d_1$ and the $d_\infty$-metric, by at least $\left(\frac{1}{n2^n}\right)^2$. A finite number $T$ of iterations are, therefore, sufficient for obtaining a solution. More specifically, $O(n)$
bisections on $\alpha$ are needed before $u_t-l_t\le \left(\frac{1}{n2^n}\right)^2$ and further improvements become theoretically impossible.}

\highlight{A pseudo-code description of this bisection approach reads as follows:\footnote{\highlight{The description focuses on finding the minimal distance $\alpha^*$. A game $v^*\in\Gamma$ with $B(v^*)=\alpha^*$ can straightforwardly be obtained from the solution to the feasibility problem $\langle \beta,\Gamma,d(\cdot), \alpha^*\rangle$.}}}

\medskip

\noindent
\highlight{\textbf{Input:} desired power index vector $\beta$, class of binary voting systems $\Gamma$, metric $d(\cdot)$\\
\textbf{Output:} minimum $d$-distance $\alpha^*$ between $\beta$ and PBI vectors induced by $\Gamma$\\[1mm]
$l_1=0$\\
$u_1=2$\\
$\alpha^*=2$\\
$\varepsilon=\left(\frac{1}{n2^n}\right)^2$\\
$t=1$\\
while $u_t-l_t>\varepsilon$\\
\hspace*{2mm} $\alpha=\frac{u_{t}-l_{t}}{2}$\\
\hspace*{2mm} $t=t+1$\\
\hspace*{2mm} solve feasibility problem $\langle \beta,\Gamma,d(\cdot), \alpha\rangle$ \\
\hspace*{2mm} if $v\in\Gamma$ such that $d(B(v)-\beta)\le \alpha$ exists\\
\hspace*{2mm} then\\
\hspace*{4mm} $u_t=d(B(v)-\beta)$, $\alpha^*=u_t$\\
\hspace*{2mm} else\\
\hspace*{4mm} $l_t=\alpha$\\
\hspace*{2mm} end if\\
end while\\
return $\alpha^*$
}

\medskip
\highlight{The \emph{feasibility problem} $\langle \beta,\Gamma,d(\cdot), \alpha\rangle$ consists of verifying whether there exists a voting system $v\in \Gamma$ such that $d(B(v),\beta)\le \alpha$.
The following ILP formulation describes it for $\Gamma=\mathcal{S}$ and the $d_1$-metric.}
Adaptations to $\mathcal{C}$ or $\mathcal{W}$ and $d_{1}'(\cdot)$ or $d_\infty(\cdot)$ involve further variables and (modified) constraints, but are otherwise very similar:

\begin{align}
  & x_S\in\{0,1\}&&\forall S\subseteq N,\label{ie_sg_start}\\
  & x_S\le x_T && \forall S\subseteq T \subseteq N,\\
  & x_\emptyset=0\\
  & x_N=1\label{ie_sg_end}\\ \noalign{\vskip5pt}
  & y_{i,S}\in\{0,1\}&&\forall 1\le i\le n,S\subseteq N\backslash\{i\},\label{ie_swing_start}\\
  & y_{i,S}=x_{S\cup\{i\}}-x_S&&\forall 1\le i\le n, S\subseteq N\backslash\{i\},\label{ie_swing_end}\\ \noalign{\vskip5pt}
  & s_i\ge 0 &&\forall 1\le i\le n,\\
  & s_i=\sum\nolimits_{S\subseteq N\backslash\{i\}} y_{i,S}&&\forall 1\le i\le n,\label{ie_count_swings}\\
  & s=\sum\nolimits_{i=1}^n s_i,&& \label{ie_sum_swings}\\\noalign{\vskip5pt}
  & \delta_i\ge 0 &&\forall 1\le i\le n,\label{ie_delta_i}\\
  & \delta_i\ge s_i- \beta_i\cdot s &&\forall 1\le i\le n,\label{ie_abs_1}\\
  & \delta_i\ge -s_i+\beta_i\cdot s &&\forall 1\le i\le n,\label{ie_abs_2}\\ \noalign{\vskip5pt}
  & \sum\nolimits_{i=1}^n \delta_i\le \alpha\cdot s. \label{ie_alpha}
\end{align}

The binary variables $x_{S}$ define a Boolean function $v$ via $\chi(S)=x_S$; inequalities (\ref{ie_sg_start})--(\ref{ie_sg_end}) ensure that they represent a simple game.
The binary auxiliary variables $y_{i,S}=x_{S\cup\{i\}}-x_S$ which are introduced in (\ref{ie_swing_start})--(\ref{ie_swing_end}) for all $i\in N$ and $\emptyset\subseteq S\subseteq N\backslash\{i\}$ satisfy $y_{i,S}=1$ if and only if coalition $S$ is a swing for voter $i$, i.e., contributes $1/2^{n-1}$ to $B_i'(v)$. \highlight{They are used in order to determine the number of swings $s_i=2^{n-1}\cdot B_i'(v)$ for each player $i$ in equation~(\ref{ie_count_swings}). 
The total number of swings $s=\sum_{i=1}^ns_i$ is defined in equation~(\ref{ie_sum_swings}). Based on this total number, the individual deviation $\delta_i = \left|s_i-\beta_i\cdot s\right|$ from the target number of swings is captured by inequalities (\ref{ie_abs_1}) and (\ref{ie_abs_2}). The feasibility of a $d_1$-distance $\alpha$ is then finally checked by introducing constraint (\ref{ie_alpha}). Namely, a simple game $v\in \mathcal{S}$ whose PBI has $d_1$-distance of $\alpha$ or less exists if and only if the feasible set defined by (\ref{ie_sg_start})--(\ref{ie_alpha}) is non-empty.}

\highlight{The answer to whether this is the case~-- and, as a by-product, some $v\in\mathcal{S}$ with distance at most $\alpha$~-- can be obtained by feeding (\ref{ie_sg_start})--(\ref{ie_alpha}) into a standard ILP software package in the required format. We have used IBM ILOG CPLEX~12.4 and the hardware described in Section~\ref{sec_results}.}

\section*{\highlight{Appendix B: Analytical PBI calculations}}
\label{sec_appendix_pbi_calculations}

\highlight{This appendix presents some technical details on the PBI computations for the sequence of weighted voting games $\{v^n\}_{n\in \mathbb{N}}$ with
$$
v^n=(q^n;w^n)=(2n-a-4;\underbrace{3,\ldots, 3}_{a \textnormal{ threes}}, \underbrace{2,\ldots, 2}_{n-a-1 \textnormal{ twos}},1),
$$
which is considered in the final paragraphs of Section~\ref{subsec:analytical}. Our first lemma determines the number of
swings in $v^n$ for each voter $i=1, \ldots, n$, that is, the cardinality of set $\{S\subseteq N\setminus \{i\}\colon {v^n}(S\cup \{i\})-{v^n}(S)=1 \}$, as a function of $a$ and $n$.}

\highlight{\begin{lemma}
  \label{lemma_swings}
The numbers of swings in ${v^n}$ are $2n-a-2$, $2n-a-4$, and $a$ for all voters with weight 3, 2, and 1, respectively.
\end{lemma}}
\begin{proof}
  \highlight{It is convenient to exploit the fact that for any $v\in \mathcal{S}$ the number of voter~$i$'s swings in $v$ and in the \emph{dual game} $v'\in \mathcal{S}$ which is obtained by setting $v'(S)=1-v(S)$ for all $S\subseteq N$ must coincide. So instead of $v^n$ consider the game ${v^n}'$ which involves identical weights but quota $q'=4$ instead of $2n-a-4$. Referring to winning and losing coalitions in ${v^n}'$ we have:}
\begin{itemize}
\item[(i)] \highlight{A voter $i$ with $w_i=3$ renders a losing coalition $S\subseteq N\setminus \{i\}$ winning by joining if either $|S|=1$ or  $S=\{j,k\}$ with $w_j=2$ and $w_k=1$. There are $n-1$ coalitions of the former and $n-a-1$ coalitions of the latter type, amounting to $2n-a-2$ swings altogether.}
\item[(ii)] \highlight{A voter $i$ with $w_i=2\subseteq N\setminus \{i\}$ renders a losing coalition $S$ winning by joining if either $S=\{w_j\}$ with $w_j=3$ or $2$, or $S=\{j,k\}$ with $w_j=2$ and $w_k=1$. There are $a+(n-a-2)$ coalitions of the former and $n-a-2$ coalitions of the latter type, amounting to $2n-a-4$ swings altogether.}
\item[(iii)] \highlight{Voter $n$ with $w_n=1$ renders a losing coalition $S\subseteq N\setminus \{i\}$ winning by joining if $S=\{j\}$ with $w_j=3$. There are $a$ such coalitions.}
\end{itemize}
\hspace{11.7cm}
\qed
\end{proof}

\highlight{Writing $\lfloor x \rfloor$ to denote the largest integer not greater than $x$, and $x \textnormal{ mod } y$ to denote the integer remainder when $x$ is divided by $y$, we have the following finding for distances in the $d_1$-metric:}

\highlight{\begin{lemma}
\label{lem:6n7specifics}
Choose
$$
a(n)=
\begin{cases}
\left\lfloor\frac{6n}{7}\right\rfloor & \textnormal{if } n \textnormal{ mod } 7 \in \{1,2,3\}, \\
\left\lfloor\frac{6n}{7}\right\rfloor-1 & \textnormal{if } n \textnormal{ mod } 7 \in \{0,4,5,6\},
\end{cases}
$$
and consider
$$
{v^n} = (2n-a-4;\underbrace{3,\ldots, 3}_{a(n) \textnormal{ threes}}, \underbrace{2,\ldots, 2}_{n-a(n)-1 \textnormal{ twos}},1).
$$
Then
$$
\lim_{n\to\infty} n\cdot\left(d_1(B({v^n}),\beta^n)\right)=\frac{1}{2}.
$$
\end{lemma}}
\begin{proof}

\highlight{Suppose that $n\textnormal{ mod } 7=0$, i.e., $n=7k+0$ for some $k\in \mathbb{N}$. Then $a=6k-1$ and
Lemma~\ref{lemma_swings} yields swing numbers of $2\cdot 7k -(6k-1) -2$, $2\cdot 7k -(6k-1) -4$, and $6k-1$ for the three voter types, respectively. This implies a total number of
$$
(6k-1)(8k-1)+(n-6k)(8k-3)+(6k-1)=56k^2-11k
$$
swings, and hence a PBI vector of
$$
B(v^{7k})=\frac{1}{56k^2-11k}\left(\underbrace{8k-1, \ldots, 8k-1}_{6k-1 \textnormal{ times}},
\underbrace{8k-3, \ldots, 8k-3}_{k \textnormal{ times}}, 6k-1  \right).
$$
This yields
$$
B_i(v^{7k})-\beta_i^{7k} =
\begin{cases}
-\frac{1}{(14k-1)k(56k-11)} & \textnormal{if }w_i=3, \\
\frac{28k-3}{(14k-1)k(56k-11)} & \textnormal{if }w_i=2, \\
-\frac{28k^2-9k+1}{(14k-1)k(56k-11)} &\textnormal{if }w_i=1,
\end{cases}
$$
and summing the absolute values of these figures up for the $6k-1$ voters with weight 3, the $k$ voters with weight~2, and the final voter $n$ one obtains
$$
\Vert B_i(v^n)-\beta^n \Vert_1=   \frac{2(28k-3)}{(14k-1)k(56k-11)}
$$
in case of $n=7k$. This number and results of the similarly tedious computations when $n\textnormal{ mod } 7=1, \ldots, 6$ are summarized in Table~\ref{table_tedious}.
For each of the seven cases one easily sees that the deviations tend to $\frac{1}{14k}$, which is equivalent to $\frac{1}{2n}$.\hfill{\qed}}

\begin{table}[t]
  \begin{center}
    \begin{tabular}{|c|c|c|c|c|}
      \hline
      {$n=$} & \multicolumn{3}{c|}{$B_i({v^n})-\beta_i$ for} & \multirow{2}{*}{$\Vert B({v^n})-\beta\Vert_1$}  \\
      $7k+$ & $w_i=3$ & $w_i=2$ &
      $w_i=1$ & \\
      \hline
0 &   $-\frac{1}{(14k-1)k(56k-11)}$ & $\frac{28k-3}{(14k-1)k(56k-11)}$ & $-\frac{28k^2-9k+1}{(14k-1)k(56k-11)}$ & $\frac{2(28k-3)}{(14k-1)k(56k-11)}$      \\ \hline
1 &     0 & $\frac{1}{2k(14k+1)}$ & $-\frac{1}{2(14k+1)}$ & $\frac{1}{14k+1}$ \\ \hline
2 &    $\frac{1}{(14k+3)(56k^2+19k+2)}$ & $\frac{7(4k+1)}{(14k+3)(56k^2+19k+2)}$ & $-\frac{28k^2+13k+1}{(14k+3)(56k^2+19k+2)}$ & $\frac{2(28k^2+13k+1)}{(14k+3)(56k^2+19k+2)}$ \\ \hline
3 &    $\frac{1}{(14k+5)(28k^2+17k+3)}$  & $\frac{2(7k+3)}{(14k+5)(28k^2+17k+3)}$ & $-\frac{2(7k^2+6k+1)}{(14k+5)(28k^2+17k+3)}$ & $\frac{4(7k^2+6k+1)}{(14k+5)(28k^2+17k+3)}$ \\ \hline
4 &    $-\frac{1}{7(2k+1)(14k^2+14k+3)}$ & $\frac{14k+5}{14(2k+1)(14k^2+14k+3)}$ & $-\frac{14k^2+7k+1}{14(2k+1)(14k^2+14k+3)}$ & $\frac{14k^2+19k+5}{7(2k+1)(14k^2+14k+3)}$ \\ \hline
5 &    $-\frac{3}{(14k+9)(56k^2+71k+21)}$ & $\frac{28k+15}{(14k+9)(56k^2+71k+21)}$ & $-\frac{28k^2+25k+6}{(14k+9)(56k^2+71k+21)}$ & $\frac{2(28k^2+43k+15)}{(14k+9)(56k^2+71k+21)}$ \\ \hline
6 &    $-\frac{1}{(14k+11)(28k^2+43k+16)}$ & $\frac{2(7k+5)}{(14k+11)(28k^2+43k+16)}$ & $-\frac{2(7k^2+9k+3)}{(14k+11)(28k^2+43k+16)}$ & $\frac{4(7k^2+12k+5)}{(14k+11)(28k^2+43k+16)}$ \\ \hline
    \end{tabular}
    \caption{$d_1$-distances between $\beta^n$ and PBI of game $v^n\in\mathcal{W}$ in Lemma~\ref{lem:6n7specifics} (with $k\in \mathbb{N}$)}
    \vspace*{-10mm}
    \label{table_tedious}
  \end{center}
\end{table}

\end{proof}

\highlight{In case of the $d_\infty$-metric, choose
$$
a(n)=\left\lfloor(n+1)/3\right\rfloor+\left\lfloor n/3\right\rfloor-1.$$
In each corresponding game $v^n$ (see Lemma~\ref{lem:6n7specifics}) roughly two thirds of the players have weight $3$ each, roughly one
third have weight $2$, and a single player has weight $1$. The games $v^n$ result in very good solutions of the inverse problem for $n<8$ and the best ones we could find for
$n\ge 8$. Using this in order to obtain an upper bound one can verify the following result in perfect analogy to Lemma~\ref{lem:6n7specifics}:}

\begin{lemma}
\highlight{The weighted voting game $v^{***}\in \mathcal{W}$ whose PBI minimizes $d_\infty$-distance to $\beta^n$ satisfies\\[-4mm]
  $$
    d_\infty\!\left(v^{***},\beta^n\right)\le b(n)=
    \begin{cases}
    \frac{8n-9}{n(4n-7)(2n-1)} & \textnormal{if } n \textnormal{ mod } 3=0,\\
    \frac{8n-23}{(4n^2-5n-8)(2n-1)} & \textnormal{if } n \textnormal{ mod } 3=1, \\
    \frac{4}{4n^2-1} & \textnormal{if } n \textnormal{ mod } 3=2
    \end{cases}
  $$\\[-4mm]
for $n\ge 8$.}
\end{lemma}
\highlight{Note that the indicated bound tends to $\frac{1}{n^2}$, i.e.,
$\lim_{n\to \infty} {b(n)}\big/{ \frac{1}{n^2}}=1.$}\\[-10mm]

\setlength{\labelsep}{-0.2cm}
{
\newcommand{\noopsort}[1]{}

\end{document}